\documentclass[11pt]{article}
\usepackage{color}
\usepackage{amssymb}
\usepackage{amsthm,array,amssymb,amscd,amsfonts,latexsym, url}
\usepackage{amsmath}
\newtheorem{theo}{Theorem}[section]

\newtheorem{prop}[theo]{Proposition}
\newtheorem{lemm}[theo]{Lemma}

\newtheorem{rema}[theo]{Remark}
\newtheorem{Defi}[theo]{Definition}

\newtheorem{question}[theo]{Question}

\newcommand{\cqfd}
{%
\mbox{}%
\nolinebreak%
\hfill%
\rule{2mm}{2mm}%
\medbreak%
\par%
}
\newfont{\gothic}{eufb10}

\date{}

\begin{document}
\title{ Abel-Jacobi map, integral Hodge classes and \\decomposition of the diagonal}
\author{Claire Voisin}
 \maketitle \setcounter{section}{-1}

\begin{flushright}  {\it \`{A} Eckart Viehweg}
\end{flushright}

\begin{abstract}
   Given a smooth projective $3$-fold $Y$, with $H^{3,0}(Y)=0$, the Abel-Jacobi map induces
   a morphism from each smooth variety parameterizing $1$-cycles in $Y$ to the intermediate Jacobian
   $J(Y)$.
   We study in this paper  the existence of families of $1$-cycles in $Y$ for which this induced morphism
   is surjective with  rationally connected general fiber, and various applications of this property.
   When $Y$ itself is rationally connected, we relate this property to the existence of an integral  homological decomposition of the diagonal.
   We also study this property for cubic threefolds, completing the work of Iliev-Markoushevich. We then
   conclude that the Hodge conjecture  holds for
   degree $4$ integral Hodge classes on fibrations into cubic threefolds over  curves,
   with restriction on singular  fibers.

\end{abstract}

\maketitle
\section{Introduction}
The following result is proved by Bloch and Srinivas \cite{blochsrinivas} as a consequence of their decomposition of the diagonal:
\begin{theo} \label{theointrotors} Let $Y$ be a  smooth projective $3$-fold  with $CH_0(Y)$ supported on a curve. Then, via the
 Abel-Jacobi map $AJ_Y$, the group of $1$-cycles homologous
to $0$ on $Y$ modulo rational equivalence is up to torsion isomorphic to the intermediate Jacobian $J(Y)$ of
$Y$,
$$J(Y):=H^3(Y,\mathbb{C})/(F^2H^3(Y)\oplus H^3(Y,\mathbb{Z})).$$
\end{theo}
This follows  from the diagonal decomposition principle
\cite{blochsrinivas} (see also \cite[II, 10.2.1]{voisinbook}). Indeed, this principle says that
if $Y$ is a smooth variety such that $CH_0(Y)$ is supported on some subvariety $W\subset Y$, there is an equality in
$CH^d(Y\times Y),\,d=dim\,Y$:
\begin{eqnarray}\label{decompintro}
N\Delta_Y=\Gamma_1+\Gamma_2,
\end{eqnarray}
where $N$ is a nonzero integer and $Z_1,\,Z_2$ are codimension $d$ cycles with
$$|Z_1|\subset D\times Y,\,D\varsubsetneq Y,\,\,\,\,|Z_2|\subset Y\times W.$$
Assuming ${\rm dim}\,W\leq 1$, the integer $N$ appearing in this decomposition is easily checked
to annihilate the kernel and cokernel of $AJ_Y$.

Note that the integer $N$ appearing above cannot be in general set to be $1$, and one
purpose of this paper is to investigate the significance of this invariant, at least if we work on the level of cycles modulo homological equivalence. We will say that $Y$ admits an
 integral homological  decomposition of the diagonal as in (\ref{decompintro}) if there is an equality as in (\ref{decompintro}) (with the same conditions on ${\rm Supp}\,Z_1$ and ${\rm Supp}\,Z_2$), with $N=1$, and
in $CH^d(Y\times Y)/{\rm hom}$.

Going further and using the theory of Bloch-Ogus
\cite{blochogus} together with  Merkurjev-Suslin theorem, Bloch and Srinivas
also prove the following:
\begin{theo} \label{blochsrinivasgriffnul} If $Y$ is a smooth projective complex variety such that
 $CH_0(Y)$ is  supported on a surface, the Griffiths group ${\rm Griff}^2(Y)=CH^2(Y)_{hom}/{\rm alg}$ is identically
$0$.
\end{theo}
The last result cannot be obtained as a consequence of the
diagonal decomposition, which only shows that under the same assumption  ${\rm Griff}^2(Y)$ is annihilated by the integer $N$ introduced above.

We finally have the following improvement of Theorem
\ref{theointrotors} :
\begin{theo} \label{introimprove}(cf. \cite{murre}) If $Y$ is a smooth projective $3$-fold such that
 $CH_0(Y)$ is  supported on a curve,
  then the
Abel-Jacobi map induces an isomorphism:
$$AJ_Y:CH_1(Y)_{hom}=CH_1(Y)_{alg}\rightarrow J(Y).$$
\end{theo}
This follows indeed from the fact that this map is surjective with  kernel of torsion
by Theorem \ref{theointrotors}, and that it can be shown by delicate arguments involving
the Merkurjev-Suslin theorem, Gersten-Quillen resolution in $K$-theory, and Bloch-Ogus theory \cite{blochogus}, that
this map induces an isomorphism on torsion.

The group on the left does not have a priori the structure of an algebraic variety, unlike the group on the right.
However it makes sense to say that $AJ_Y$ is algebraic, meaning that for any smooth
algebraic variety $B$, and any codimension $2$-cycle
$Z\subset B\times Y$, with $Z_b\in CH^2(Y)_{hom}$ for any $b\in B$, the induced map
$$\phi_Z:B\rightarrow J(Y),\,b\mapsto AJ_Y(Z_b),$$
is a morphism of algebraic varieties.

Theorem \ref{introimprove}  may
give the feeling that for such threefolds,  the study of $1$-cycles  parallels that of $0$-cycles on curves.
Still the structure of
the Abel-Jacobi map on families of $1$-cycles on such threefolds is mysterious, in contrast to what happens in the curve case, where Abel's theorem shows that the Abel-Jacobi
map on the family of effective  $0$-cycles of large degree has fibers isomorphic to projective spaces. Another goal of  this paper is to underline substantial  differences between $1$-cycles on threefolds with small $CH_0$ on one side  and $0$-cycles on curves on the other side.

There are for example two natural questions left open by the above results :

\begin{question} \label{question1} Let $Y$ be a smooth projective threefold, such that
$AJ_Y:CH_1(Y)_{alg}\rightarrow J(Y)$ is surjective.
Is there a family of codimension $2$-cycles homologous to
$0$ on $Y$, parameterized by $J(Y)$, say $Z\subset J(Y)\times Y$,  such that the induced morphism
$\phi_Z:J(Y)\rightarrow J(Y)$
is the identity?
\end{question}

Note that the surjectivity assumption is conjecturally implied by the vanishing $H^3(Y,\mathcal{O}_Y)=0$, via
the generalized Hodge conjecture (cf. \cite{grothendieck}).

There are two other related questions that we wish to partially investigate  in this paper:
\begin{question} \label{question2} Is the following property (*) satisfied by $Y$?

(*) There exists a family of codimension $2$-cycles homologous to
$0$ on $Y$, parameterized by a smooth projective variety
$B$, say $Z\subset B\times Y$, with $Z_b\in CH^2(Y)_{hom}$ for any $ b\in B$, such that the induced morphism
$\phi_Z:B\rightarrow J(Y)$
is surjective with rationally connected general fiber.
\end{question}
This question has been solved by Iliev-Markoushevich (\cite{ilievmarku}, see also \cite{HRSabeljacobi} for
similar results obtained independently)
in the case where $Y$ is a smooth cubic threefold in $\mathbb{P}^4$.
Obviously a positive answer to Question \ref{question1} implies a positive answer to Question \ref{question2}, as we can then just take
$B=J(X)$.
However, it seems that Question \ref{question2} is more natural, especially if we go to the following stronger version :
Here we choose an integral cohomology class $\alpha\in H^4(Y,\mathbb{Z})$. Assuming $CH_0(Y)$
is supported on a curve, one has $H^2(Y,\mathcal{O}_Y)=0$ by Mumford's theorem \cite{mumford}
(see also \cite[II,10.2.2]{voisinbook}) and thus $\alpha$ is a Hodge class.
Introduce the torsor $J(Y)_\alpha$  defined as follows:
the Deligne cohomology group $H^4_D(Y,\mathbb{Z}(2))$ is an extension
$$0\rightarrow J(Y)\rightarrow H^4_D(Y,\mathbb{Z}(2)) \stackrel{o}{\rightarrow} Hdg^4(Y,\mathbb{Z})\rightarrow0,$$
where $o$ is the natural map from Deligne to Betti cohomology.
Define \begin{eqnarray}\label{eqn7avrildeligne}J(Y)_\alpha:=o^{-1}(\alpha).
 \end{eqnarray}
 By definition, the Deligne cycle class (cf. \cite{EV}, \cite[I,12.3.3]{voisinbook}), restricted to codimension $2$ cycles of
class $\alpha$,  takes value in $J(Y)_\alpha$. For any family
 of $1$-cycles ${Z}\subset B\times Y$   of class $[Z_t]=\alpha,\,t\in B$ on $Y$
 parameterized by an algebraic variety
  $B$, the Abel-Jacobi map (or rather the Deligne cycle class map) of $Y$ induces a morphism
 complex algebraic varieties:
 $$\phi_{Z}:B\rightarrow J(Y)_\alpha.$$
\begin{question}\label{Q3} Is the following property (**) satisfied by $Y$?

(**) For any degree $4$ integral cohomology class $\alpha$ on $Y$,  there is a {\it canonically defined} variety $B_\alpha$,
together with a family of  $1$-cycles of class
$\alpha$ on $Y$, parameterized by
$B_\alpha$, say $Z_\alpha\subset B_\alpha\times Y$, with $[Z_{\alpha,b}]=\alpha$ in $ H^4(Y,\mathbb{Z})$ for any $ b\in B$, such that the morphism
$\phi_{Z_\alpha}:B_\alpha\rightarrow J(Y)_\alpha$ induced  by Deligne cycle class
is surjective with rationally connected general fiber.
\end{question}
Let us comment on the importance of Question \ref{Q3} in relation with the Hodge conjecture with integral coefficients for degree $4$ Hodge classes, (which has been recently related in \cite{ctvoisin} to
degree $3$ unramified cohomology with finite coefficients): The important point here is that we
want to consider fourfolds fibered over curves, or families of  threefolds parameterized by a curve. Property (**) essentially
says that
property (*), being satisfied in a canonical way, holds in family. When we work
 in families,  the necessity  to look at all twistors
$J(Y)_\alpha$, and not only at $J(Y)$ becomes obvious: for fixed $Y$ the twisted Jacobians are all
isomorphic (maybe not canonically); this is not true in families: non trivial twistors
$\mathcal{J}_\alpha$  can appear. We will give  more precise  explanations in subsection \ref{subsectionintro}
of this introduction and explain one application of this property to the Hodge conjecture for degree $4$ integral Hodge classes on fourfolds fibered over curves.

Our results in this paper are of two kinds. First of all, we extend the results of
\cite{ilievmarku} and answer affirmatively Question \ref{Q3} for cubic threefolds. As a consequence, we get the following
result (Theorem \ref{cubique}):
\begin{theo} \label{cubiqueintro}Let  $f:X\rightarrow \Gamma$ be a  fibration over a curve with general
fiber a smooth cubic threefold. If the fibers
 of $f$ have at worst  ordinary quadratic  singularities,  then the Hodge conjecture holds for
 degree $4$ integral Hodge classes on $X$.

\end{theo}
By the results of \cite{ctvoisin}, this result can also be translated into a statement concerning
degree $3$ unramified cohomology with $\mathbb{Q}/\mathbb{Z}$-coefficients of such fourfolds (see section
\ref{subsectionintro}).

Our second result relates Question \ref{question2} to the existence of  an homological integral decomposition of the diagonal as in (\ref{decompintro}).
\begin{theo} \label{theoechant} i) Let $Y$ be a smooth projective  $3$-fold.
If $Y$ admits  an integral homological  decomposition of the diagonal as in (\ref{decompintro}),  then  $H^4(Y,\mathbb{Z})$ is generated by classes of algebraic cycles,
$H^*(Y,\mathbb{Z})$ has no torsion and  $Y$ satisfies condition (*).

ii) Conversely,  assume  $CH_0(Y)=\mathbb{Z}$, $H^4(Y,\mathbb{Z})$ is generated by classes of algebraic cycles and $H^*(Y,\mathbb{Z})$ is without torsion (these conditions are satisfied by a rationally connected threefold with no torsion in $H^3(Y,\mathbb{Z})$).

If furthermore  the intermediate Jacobian of
$Y$ admits a $1$-cycle $\Gamma$ of class $\frac{[\Theta]^{g-1}}{(g-1)!}$, $g=dim\,J(Y)$,
then
  condition (*)
implies
that $Y$ admis  an integral homological  decomposition of the diagonal as in (\ref{decompintro}).
\end{theo}

Here we recall that the intermediate Jacobian $J(Y)$ is naturally a principally polarized abelian variety, the polarization $\Theta$ being given by the intersection form on $H^3(Y,\mathbb{Z})\cong H_1(J(Y),\mathbb{Z})$.

We will also prove the following relation between
Question \ref{question1} and \ref{question2} (cf. Theorem \ref{theonouveau}):
\begin{theo} Assume that Question \ref{question2} has a positive answer for $Y$ and that
the intermediate Jacobian of
$Y$ admits a  $1$-cycle $\Gamma$ of class $\frac{[\Theta]^{g-1}}{(g-1)!}$, $g=dim\,J(Y)$. Then
Question \ref{question1} also has an affirmative answer for $Y$.
\end{theo}

The paper is organized as follows: in section \ref{secfibrecub}, we give a positive
answer to Question \ref{question2} for very general cubic threefolds. We deduce from this Theorem
\ref{cubiqueintro}.
Section \ref{vraiesecdecomp} is devoted to various relations between Question \ref{question1} and \ref{question2}
and the relation between these questions and the homological decomposition of the diagonal with integral coefficients, in the spirit of Theorem \ref{theoechant}.

{\bf Thanks.} I wish to thank Jean-Louis Colliot-Th\'{e}l\`{e}ne for many discussions leading to the
present form of the paper and for numerous questions and suggestions. In particular, Theorem \ref{cubiqueintro} answers
one of his questions. I also thank Dimitri Markoushevich and Jason Starr for helpful communications.

\section{\label{subsectionintro} Preliminaries on integral Hodge classes and unramified cohomology}

We give in this section  a description of results and notations from the
earlier papers \cite{voisinuniruled}, \cite{ctvoisin} which will be used later on in the paper.
Let $X$ be a smooth projective complex variety. We denote  $Z^{2i}(X)$ the quotient
$Hdg^{2i}(X,\mathbb{Z})/H^{2i}(X,\mathbb{Z})_{alg}$ of the group  of degree $2i$ integral Hodge classes
on $X$ by the subgroup consisting of cycle classes.
The group $Z^{2i}(X)$ is a birational invariant of $X$ for $i=2,\,d-1$, where $d=dim\,X$
(cf. \cite{soulevoisin}). It is of course trivial in degrees
 $2i=0,\,2d,\,2$, where the Hodge conjecture holds for integral Hodge classes.

 The paper \cite{ctvoisin} focuses on $Z^4(X)$ and relates its torsion to degree $3$ unramified cohomology
 of $X$ with finite coefficients.
 Recall  that for any constant
  coefficients sheaf $A$ on $X$, unramified cohomology $H^i_{nr}(X,A)$ with coefficients in $A$ was introduced
   in  \cite{colliotojanguren} : denote by $X_{Zar}$
 the variety $X$ endowed with the Zariski topology, while $X$ will be considered as endowed with the classical topology. The identity map
 $$\pi:X\rightarrow X_{Zar}$$
 is continuous and allows Bloch and Ogus \cite{blochogus} to introduce sheaves
 $\mathcal{H}^l(A)$
 on $X_{Zar}$ defined by $$\mathcal{H}^l(A):=R^l\pi_*A.$$
 \begin{Defi} (Ojanguren-Colliot-Th\'{e}l\`{e}ne)  Unramified cohomology $H^i_{nr}(X,A)$ of $X$ with coefficients in $A$ is  the group $H^0(X_{Zar},\mathcal{H}^i(A))$.
 \end{Defi}
 This is a birational invariant of $X$, as shown by the Gersten-Quillen resolution for the sheaves
 $\mathcal{H}^i$  proved by Bloch-Ogus \cite{blochogus}. We refer to \cite{ctvoisin} for the description of
 other birational invariants constructed from the cohomology of the sheaves $\mathcal{H}^i$.

 The following result is proved in \cite{ctvoisin} (assuming the Bloch-Kato conjecture announced by Voevodsky,
 or restricting to $n=2$ for an unconditional result):
 \begin{theo} \label{ctVtheo}There is an exact sequence for any $n$:
 $$0\rightarrow H^3_{nr}(X,\mathbb{Z})\otimes \mathbb{Z}/n\mathbb{Z}\rightarrow H^3_{nr}(X,\mathbb{Z}/n\mathbb{Z})\rightarrow Z^4(X)[n]\rightarrow 0,$$
 where the notation $[n]$ means that we take the $n$-torsion. There is an exact sequence:
 $$0\rightarrow H^3_{nr}(X,\mathbb{Z})\otimes \mathbb{Q}/\mathbb{Z}\rightarrow H^3_{nr}(X,\mathbb{Q}/\mathbb{Z})\rightarrow Z^4(X)_{tors}\rightarrow 0,$$
 where the subscript ``tors'' means that we consider the torsion part.
 Furthermore,  the first term   is $0$ if
 $CH_0(X)$ is supported on a surface.
 \end{theo}

 Concerning the vanishing of the group $Z^4(X)$, the following is proved in
 \cite{voisinuniruled}. It will be used in section \ref{secdecomp}
\begin{theo}\label{theovoisinuniruled} Let $Y$ be a smooth projective threefold which is either uniruled
or Calabi-Yau. Then $Z^4(Y)=0$, that is, any integral degree $4$ Hodge class on $Y$ is algebraic. In particular, if $Y$ is uniruled with
$H^2(Y,\mathcal{O}_Y)=0$, any integral degree $4$ cohomology class on $Y$ is algebraic.
\end{theo}

 We will now be interested in the group
  $Z^4(X)$, where  $X$ is a smooth projective
  $4$-fold, $f:X\rightarrow \Gamma$ is  a   surjective morphism to a smooth curve
  $\Gamma$, whose general fiber
 $X_t$  satisfies
 $H^3(X_t,\mathcal{O}_{X_t})=H^2(X_t,\mathcal{O}_{X_t})=0$.
  For  $X_t,\,t\in \Gamma,$ as above, the intermediate Jacobian
  $J(X_t)$ is an abelian variety, as a consequence of the vanishing  $H^3(X_t,\mathcal{O}_{X_t})=0$ (cf. \cite[I,12.2.2]{voisinbook}).
  For any class
 $$\alpha\in H^4(X_t(\mathbb{C}),\mathbb{Z})=Hdg^4(X_t,\mathbb{Z}),$$ we introduced above a torsor
 $J(X_t)_\alpha$ under  $J^2(X_t)$ which is an algebraic variety non canonically isomorphic to
 $J(X_t)$.

 This construction works in family on the Zariski open set $\Gamma_0\subset \Gamma$, over which $f$
 is smooth. There is thus  a family of abelian varieties
  $\mathcal{J}\rightarrow \Gamma_0$, and for any locally constant section  $\alpha $
 of the locally constant system  $R^4f_*\mathbb{Z}$, we get the twisted
  family $\mathcal{J}_\alpha\rightarrow \Gamma_0$, defined by the obvious extension of
  formula (\ref{eqn7avrildeligne}) in the relative setting. The construction of these families in the analytic setting is easy. The fact that the resulting families are algebraic is much harder.

 Given  a smooth algebraic variety  $B$, a morphism $g:B\rightarrow \Gamma$ and
 a  family of relative  $1$-cycles  $\mathcal{Z}\subset B\times_\Gamma X$ of class
 $[Z_t]=\alpha_t\in H^4(X_t,\mathbb{Z})$, the  relative Abel-Jacobi map (or rather Deligne cycle class)
 gives a  morphism
 $\phi_\mathcal{Z}:B_0\rightarrow \mathcal{J}_\alpha$
 over  $\Gamma_0$, where  $B_0:=g^{-1}(\Gamma_0)$.

The following result, which illustrates the importance of condition
 (**) as opposed to condition (*), appears  in \cite{ctvoisin}. We recall the proof here, as we will need it to prove  a slight improvement of the criterion, which will be used  in section \ref{seccub}. Here we assume that ${X}$ is a smooth projective $4$-fold, and that $f:X\rightarrow \Gamma$ is a morphism to a smooth curve.
\begin{theo}\label{theocritereabeljacobi}  Assume  $f:{X}\rightarrow \Gamma$ satisfies the following assumptions:
 \begin{enumerate}
 \item\label{itemref0}  The smooth fibers ${X}_t$ of $f$
   satisfy the vanishing conditions
 $H^3({X}_t,\mathcal{O}_{{X}_t})=H^2({X}_t,\mathcal{O}_{{X}_t})=0$.
 \item \label{itemref1} The smooth fibers ${X}_t$ have
  no torsion in $H^3({X}_t(\mathbb{C}),\mathbb{Z})$.
     \item \label{itemref2} The singular fibers of  $f$ are reduced with at worst
     ordinary quadratic  singularities.
     \item \label{itemref3} For any section  $\alpha$ of $R^4f_*\mathbb{Z}$, there exists a variety $g_\alpha:B_\alpha\rightarrow \Gamma$  and a family of relative $1$-cycles  $\mathcal{Z}_\alpha\subset B_\alpha\times_\Gamma X$ of class
 $\alpha$, such that the  morphism
 $\phi_{\mathcal{Z}_\alpha}:B\rightarrow \mathcal{J}_\alpha$
     is surjective     with rationally connected general fiber.
     \end{enumerate}
     Then the Hodge conjecture is true for integral Hodge classes of degree $4$ on
     $\mathcal{X}$.
\end{theo}
{\bf Proof.}
 Following Zucker \cite{zucker}, an integral  Hodge class $\tilde{\alpha}\in Hdg^4(X,\mathbb{Z})\subset H^4(X,\mathbb{Z})$ induces a section $\alpha$ of the constant system $R^4f_*\mathbb{Z}$ which admits a canonical lift
 in the family of twisted Jacobians $\mathcal{J}_\alpha$. This lift is an {\it algebraic} section   $\sigma:\Gamma\rightarrow\mathcal{J}_\alpha$
of the structural map  $\mathcal{J}_\alpha \rightarrow \Gamma$.

 Recall that we have by  hypothesis the morphism
 $$\phi_{\mathcal{Z}_\alpha}:B_\alpha\rightarrow \mathcal{J}_\alpha$$
 which is algebraic,
      surjective    with rationally connected general fiber. We can now replace $\sigma$ by a
      $1$-cycle $\Sigma=\sum_in_i\Sigma_i$ rationally equivalent to it in $\mathcal{J}_\alpha$, in such a way
      that the fibers of $\phi_{\mathcal{Z}_\alpha}$ are rationally connected over the general points of
      each component $\Sigma_i$ of ${\rm Supp}\,\Sigma$.

      According to  \cite{GHS}, the morphism $\phi_{\mathcal{Z}_\alpha}$
      admits a lifting over each $\Sigma_i$, which
      provides curves $\Sigma'_i\subset B_\alpha$.

       Recall next that $B_\alpha$ parameterizes  a family
       $\mathcal{Z}_\alpha\subset B_\alpha\times_\Gamma X$ of relative $1$-cycles of class
         $\alpha$ in the fibers $X_t$. We can restrict  this family to each
         $\Sigma'_i$, getting $\mathcal{Z}_{\alpha,i}\subset \Sigma'_i\times_\Gamma X$.
       Consider the $1$-cycle
       $$Z:=\sum_in_ip_{i*}\mathcal{Z}_{\alpha,i}\subset\Gamma\times_\Gamma X=X,$$
       where $p_i$ is the restriction to $\Sigma'_i$ of $p:B_\alpha\rightarrow \Gamma$.
       Recalling that $\Sigma$ is rationally equivalent to $\sigma(\Gamma)$ in
       $\mathcal{J}_\alpha$, we find that
      the normal function  $\nu_Z$ associated to   $Z$, defined by
      $$\nu_Z(t)=AJ_{X_t}(Z_{\mid X_t})$$
       is equal to  $\sigma$.
       We then deduce from \cite{griffiths} (see also \cite[II,8.2.2]{voisinbook}), using
        the Leray spectral sequence of $f_U:X_U\rightarrow U$ and hypothesis  \ref{itemref1}, that the cohomology classes $[Z]\in H^4(X,\mathbb{Z})$ of $Z$ and
      $\tilde{\alpha}$ coincide on any open set of the form $X_U$, where  $U\subset \Gamma_0$ is an affine
      open set over which
      $f$ is smooth.

       On the other hand, the kernel of the restriction map
       $H^4(X,\mathbb{Z})\rightarrow H^4(U,\mathbb{Z})$ is generated
      by the  groups $i_{t*}H_4(X_t,\mathbb{Z})$ where  $t\in \Gamma\setminus U$, and
      $i_t:X_t\rightarrow X$ is the inclusion map. We conclude using assumption  \ref{itemref2}
      and the fact that the general fiber of  $f$ has $H^2(X_t,\mathcal{O}_{X_t})=0$, which
      imply that all fibers  $X_t$ (singular or not)  have their degree $4$ integral homology generated by
      homology classes of algebraic cycles; indeed, it follows from this and the previous conclusion
      that   $[Z]-\tilde{\alpha}$ is algebraic,
      so that $\tilde{\alpha}$ is also algebraic.

\cqfd

\begin{rema}{\rm It is also possible in this proof, instead of moving the section $\sigma$ to
a $1$-cycle in general position, to use Theorem \ref{hoxu} below, which also guarantees that in fact $\sigma$ itself lifts
to $B_\alpha$.}
\end{rema}
\begin{rema} {\rm By Theorem \ref{theovoisinuniruled}, if $X_t$ is uniruled and $H^2(X_t,\mathcal{O}_{X_t})=0$  (a strengthening of assumption \ref{itemref0}), then any degree $4$ integral cohomology class $\alpha_t$
on $X_t$
is algebraic on $X_t$. Together with
 Bloch-Srinivas results \cite{blochsrinivas} on the surjectivity of the Abel-Jacobi map
 for codimension $2$-cycles under these assumptions, this shows  that pairs $(B_\alpha,Z_\alpha)$
 with surjective $\phi_{\mathcal{Z}_\alpha}:B_\alpha\rightarrow \mathcal{J}_\alpha$
     exist. In this case, the strong statement in assumption \ref{itemref3}  is thus the rational connectedness of the fibers.}
\end{rema}
In order to study cubic threefolds fibrations, we will need a technical strenghtening of Theorem \ref{theocritereabeljacobi}.
We start with a smooth projective morphism $f:\mathcal{X}\rightarrow T$ of relative dimension $3$.
We assume as before that smooth fibers $\mathcal{X}_t$ have no torsion in
$H^3(\mathcal{X}_t,\mathbb{Z})$ and have $H^3(\mathcal{X}_t,\mathcal{O}_{\mathcal{X}_t})=H^2(\mathcal{X}_t,\mathcal{O}_{\mathcal{X}_t})=0$. As before, for any global section $\alpha$ of $R^4f_*\mathbb{Z}$, we have the twisted family of intermediate Jacobians
$\mathcal{J}_\alpha\rightarrow T$, where $T$ is smooth quasi-projective.
\begin{theo}\label{variantcritere} Assume that:

(i)  The local system $R^4f_*\mathbb{Z}$
is trivial.

 (ii) For any section  $\alpha$ of $R^4f_*\mathbb{Z}$, there exists a variety $g_\alpha:B_\alpha\rightarrow T$  and a family of relative $1$-cycles  $\mathcal{Z}_\alpha\subset B_\alpha\times_T \mathcal{X}$ of class
 $\alpha$, such that the  morphism
 $\phi_{\mathcal{Z}_\alpha}:B_\alpha\rightarrow \mathcal{J}_\alpha$
     is surjective     with rationally connected general fibers.

     Then for any smooth projective curve $\Gamma$, any morphism
     $\phi:\Gamma_0\rightarrow T$, where $\Gamma_0$ is a Zariski open set of $\Gamma$, and any smooth projective model $\psi:X\rightarrow \Gamma$ of $\mathcal{X}\times_T\Gamma_0$,  if furthermore $g$ has only nodal fibers, one has $Z^4(X)=0$.
\end{theo}
{\bf Proof.} We mimick the proof of the previous theorem. We thus only have to prove that for
any Hodge class $\tilde{\alpha}$ on $X$,
 inducing   by restriction to the fibers of $\psi$ a section $\alpha$ of $R^4\psi_*\mathbb{Z}$, the section
 $\sigma_{\tilde{\alpha}}$ of $\mathcal{J}_{\Gamma,\alpha}\rightarrow \Gamma_0$ induced by $\tilde{\alpha}$ lifts to a family of $1$-cycles of class $\alpha$ in the fibers of $\psi$.
 Observe that by triviality of the local system $R^4f_*\mathbb{Z}$ on $T$, the section $\alpha$ extends to
 a section of $R^4f_*\mathbb{Z}$ on $T$. We then have by assumption the family of $1$-cycles
 $\mathcal{Z}_\alpha\subset B_\alpha\times_T \mathcal{X}$ parameterized by $B_\alpha$, with the property that the general fiber of the induced Abel-Jacobi map $\phi_{\mathcal{Z}_\alpha}:B_\alpha\rightarrow \mathcal{J}_\alpha$ is rationally connected. As $\phi^*\mathcal{J}_\alpha=\mathcal{J}_{\Gamma,\alpha}$,
 we can see the pair $(\phi,\sigma_{\tilde{\alpha}})$ as a general morphism from $\Gamma_0$ to $\mathcal{J}_\alpha$. The desired family of $1$-cycles follows then from the existence of a
 lift of $\sigma_{\tilde{\alpha}}$ to $B_\alpha$, given by the following Theorem \ref{hoxu}, due to
 Hogadi and Xu.

\cqfd
\begin{theo} \label{hoxu} Let  $f:Z\rightarrow W$ be a surjective projective morphism between
 smooth varieties over  $\mathbb{C}$.  Assume the general fiber of $f$ is rationally connected. Then, for any rational map
$g:C\dashrightarrow W$, where $C$ is a smooth curve, there exists  a rational lift $\tilde{g}:C\rightarrow Z$
of  $g$ in $Z$.
\end{theo}

When the image curve $g(C)$ is in  general position, so that the fiber of $f$ over the general point of
$Im\,g$ is actually rationally connected, this result is a  consequence of \cite{GHS}. The degenerate case
is a
 consequence of the work of  Hogadi
and Xu \cite{hoxu} combined with \cite{GHS}. Indeed, they prove that the variety $C_0\times _WZ$ contains a subvariety fibered over $C_0$ with  rationally connected general fibers, where $C_0$ is the definition locus of $g$.
\cqfd
\section{\label{secfibrecub}On the fibres of the Abel-Jacobi map for the cubic threefold}
The papers   \cite{tikhomarku}, \cite{ilievmarku}, \cite{HRSabeljacobi} are devoted to the study
of the  morphism induced by the
Abel-Jacobi map,
from the  family of curves of small degree in a cubic threefold in $\mathbb{P}^4$
to its intermediate Jacobian.
  In degree $4$, genus $0$, and degree  $5$, genus $0$ or $1$, it is proved that these
   morphisms are  surjective with rationally connected fibers, but it is known that this is not true
      in degree $3$ (and any genus), and, to our knowlege,  the case of degree   $6$ has not been studied.
As is clear from the proof of Theorem \ref{theocritereabeljacobi}, we need for the proof of Theorem
\ref{cubique}  to have such a statement for at least one
canonically defined family of curves  of degree divisible by  $3$.
The following result provides such a statement.
\begin{theo}\label{propo61ratcon} Let  $Y\subset \mathbb{P}^4_{\mathbb{C}}$ be a general cubic  hypersurface. Then the map  $\phi_{6,1}$
induced  by the Abel-Jacobi map  $AJ_Y$, from  the family $M_{6,1}$ of elliptic curves of
 degree $6$ contained in  $Y$ to $J(Y)$, is surjective with rationally connected general fiber.
\end{theo}
\begin{rema}\label{rema14fevrier}{\rm What we call here  the family of elliptic curves is a desingularization of the closure, in the Hilbert scheme of  $Y$, of the family
parameterizing smooth  degree $6$ elliptic curves which are nondegenerate in
 $\mathbb{P}^4$.}
\end{rema}
\begin{rema}{\rm As $J(Y)$ is not uniruled, an equivalent formulation of the result is the following (cf. \cite{HRSabeljacobi}): the map $\phi_{6,1}$ is dominating and identifies birationally
 to the maximal rationally connected fibration (see \cite{komimo}) of $M_{6,1}$.}
 \end{rema}
{\bf Proof of Theorem  \ref{propo61ratcon}.}  Notice that it   suffices to prove the
statement  for very general  $Y$, which we will assume  occasionally.
One easily checks
that  $M_{6,1}$ is for general  $Y$  irreducible  of dimension $12$.  Let
$E\subset Y$ be a general nondegenerate smooth  degree $6$ elliptic curve.
Then there exists a smooth  $K3$ surface  $S\subset Y$, which is a member of the linear system $|
\mathcal{O}_Y(2)|$, containing  $E$. The line bundle
$\mathcal{O}_S(E)$ is then generated by two  sections. Let
$V:=H^0(S,\mathcal{O}_S(E))$ and consider the rank $2$ vector bundle
$\mathcal{E}$ on  $Y$ defined by  $\mathcal{E}=\mathcal{F}^*$, where
$\mathcal{F}$ is the kernel of the   surjective evaluation map
\begin{eqnarray}
\label{exactepremiermars}V\otimes\mathcal{O}_Y\rightarrow
\mathcal{O}_S(E). \end{eqnarray}
 One verifies by dualizing the
exact sequence  (\ref{exactepremiermars}) that $H^0(Y,\mathcal{E})$ is of
dimension $4$, and that  $E$ is recovered as the  zero locus of a
 transverse section of  $\mathcal{E}$. Furthermore, the bundle $\mathcal{E}$ is stable
  because $E$ is nondegenerate. One concludes that,  denoting   $M_9$ the moduli space  (of dimension $9$)
of stable rank $2$ bundles  on
$Y$ with Chern classses  $$c_1(\mathcal{E})=\mathcal{O}_Y(2),\,\,deg\,c_2(\mathcal{E})=6,$$
one has a  dominating rational map
$$\phi:M_{6,1}\dashrightarrow M_9$$
whose general fiber is isomorphic to  $\mathbb{P}^3$ (more precisely, the fiber over
$[\mathcal{E}]$ is the projective space $\mathbb{P}(H^0(Y,\mathcal{E}))$). It follows that
the  maximal rationally connected (MRC) fibration  of   $M_{6,1}$  (cf. \cite{komimo})
factorizes through $M_9$.

Let us consider the subvariety  $D_{3,3}\subset M_{6,1}$ parameterizing the singular elliptic curves
of degree  $6$ consisting of two rational  components of degree $3$ meeting  in two points.
The family  $M_{3,0}$ parameterizing rational curves of degree  $3$ is birationally
 a $\mathbb{P} ^2$-fibration over the Theta divisor  $\Theta\subset J(Y)$ (cf. \cite{HRS}).
More precisely,  a general rational curve $C\subset Y$  of degree  $3$ is contained
in a smooth  hyperplane section  $S\subset Y$, and deforms in a linear system
of dimension $2$ in  $S$. The map  $\phi_{3,0}:M_{3,0}\rightarrow J(Y)$
induced  by the Abel-Jacobi map of  $Y$  has its image equal to  $\Theta\subset J^2(Y)$
and its fiber passing through  $[C]$ is the  linear system $\mathbb{P}^2$ introduced above. This allows
to describe  $D_{3,3}$ in the following way: the Abel-Jacobi map
of $Y$, applied to  each  component $C_1,\,C_2$ of a general curve  $C_1\cup C_2$  parameterized by $D_{3,3}$ takes value in the symmetric product
$S^2\Theta$ and its fiber passing through  $[C_1\cup C_2]$ is described
by the choice of two curves in the  corresponding linear systems $\mathbb{P}^2_1$ and $\mathbb{P}^2_2$. These two curves
have to meet  in two points of the elliptic curve $E_3$ defined as the complete intersection
of the two cubic surfaces $S_1$ and $S_2$. This fiber  is thus  birationally  equivalent
to  $S^2E_3$.

 Observe that the existence of this subvariety  $D_{3,3}\subset M_{6,1}$ implies the surjectivity of $\phi_{6,1}$ because
the  restriction of  $\phi_{6,1}$  to  $D_{3,3}$ is the composition of the above-defined surjective  map
$$\chi: D_{3,3}\rightarrow S^2\Theta$$
and  of the sum map $S^2\Theta\rightarrow J(Y)$. Obviously the sum map is  surjective.

We will show more precisely:
\begin{lemm}\label{lem8mars} The map  $\phi$ introduced above is generically defined  along
$D_{3,3}$ and  $\phi_{\mid D_{3,3}}:D_{3,3}\dashrightarrow M_9$ is dominating. In particular $D_{3,3}$ dominates the base of the MRC fibration of $M_{6,1}$.
\end{lemm}
To construct other  rational curves in  $M_{6,1}$, we now  make the following observation: if $E$ is a nondegenerate degree $6$ elliptic curve in  $\mathbb{P}^4$,
let us choose  two distinct line bundles  $l_1$, $l_2$ of degree  $2$ on  $E$ such that
$\mathcal{O}_E(1)=2l_1+l_2$ in  ${\rm Pic}\,E$. Then $(l_1,l_2)$ provides an embedding of
$E$ in  $\mathbb{P}^1\times \mathbb{P}^1$, and denoting $h_i:=pr_i^*\mathcal{O}_{\mathbb{P}^1}(1)\in {\rm Pic}\,(\mathbb{P}^1\times \mathbb{P}^1)$, the line bundle  $l=2h_1+h_2$ on $\mathbb{P}^1\times \mathbb{P}^1$
has the property that
$l_{\mid E}=2l_1+l_2=\mathcal{O}_E(1)$ and the  restriction map
$$H^0(\mathbb{P}^1\times \mathbb{P}^1,l)\rightarrow H^0(E,\mathcal{O}_E(1))$$
is an  isomorphism. The original morphism from  $E$ to  $\mathbb{P}^4$ is given
by a base-point free hyperplane in  $H^0(E,\mathcal{O}_E(1))$. For a generic choice of
$l_2$, this  hyperplane also provides a base-point free  hyperplane in  $H^0(\mathbb{P}^1\times \mathbb{P}^1,l)$,
hence a  morphism $\phi:\mathbb{P}^1\times \mathbb{P}^1\rightarrow \mathbb{P}^4$, whose  image is a  surface
$\Sigma$ of degree $4$.
The residual curve of  $E$  in the intersection $\Sigma\cap Y$ is a curve
in the linear system  $|4h_1+h_2|$ on  $\mathbb{P}^1\times \mathbb{P}^1$. This is thus a curve  of
degree  $6$ in
$\mathbb{P}^4$ and of genus  $0$.

Conversely, if we start from a smooth general genus $0$ and degree $6$ curve $C$ in
$\mathbb{P}^4$, we claim that there exists a  $\mathbb{P}^1$ parameterizing  morphisms
$\phi:\mathbb{P}^1\times \mathbb{P}^1\rightarrow \mathbb{P}^4$ as above, such that
$C\subset \phi(\mathbb{P}^1\times \mathbb{P}^1)$ is the image by $\phi$ of a curve in
the linear system  $| 4h_1+h_2|$ on $\mathbb{P}^1\times \mathbb{P}^1$. To see this, we note
that the morphism from  $\mathbb{P}^1$ to $\mathbb{P}^4$
 given up to the action of  $Aut\,\mathbb{P}^1$ by  $C\subset \mathbb{P}^4$ provides a
 base-point free linear subsystem $W\subset | \mathcal{O}_{\mathbb{P}^1}(6)|$
of rank  $5$. Let us choose a  hyperplane
$H\subset | \mathcal{O}_{\mathbb{P}^1}(6)|$ containing $W$. $W$ being given, such hyperplanes $H$ are
parameterized by a  $\mathbb{P}^1$. $H$ being chosen generically, there exists a  unique
embeding $C\rightarrow \mathbb{P}^1\times \mathbb{P}^1$ as a curve in the  linear system
$|4h_1+h_2|$, such that
$H$ is recovered as the image of the  restriction map
$$H^0(\mathbb{P}^1\times \mathbb{P}^1,2h_1+h_2)\rightarrow H^0(\mathbb{P}^1,\mathcal{O}_{\mathbb{P}^1}(6)).$$
To see this last fact, notice that the embedding of  $\mathbb{P}^1$ into $\mathbb{P}^1\times \mathbb{P}^1$
as a curve in the linear system $| 4h_1+h_2|$ is determined up to the action of $Aut\,\mathbb{P}^1$ by the choice of a
degree $4$ morphism from  $\mathbb{P}^1$ to  $\mathbb{P}^1$, which is equivalent to the data of
a rank $2$ base-point free linear system  $W'\subset H^0(\mathbb{P}^1,\mathcal{O}_{\mathbb{P}^1}(4))$.
The condition that $H$  is equal to the image of the restriction map above is equivalent to the fact
 that $$H=W'\cdot H^0(\mathbb{P}^1,\mathcal{O}_{\mathbb{P}^1}(2)):=Im\,(W'\otimes H^0(\mathbb{P}^1,\mathcal{O}_{\mathbb{P}^1}(2))\rightarrow
H^0(\mathbb{P}^1,\mathcal{O}_{\mathbb{P}^1}(6)),$$
where the map is given by multiplication of sections.
Given $H$, let us set $$W':=\cap_{t\in H^0(\mathbb{P}^1,\mathcal{O}_{\mathbb{P}^1}(2))}[H:t],$$ where
$[H:t]:=\{w\in H^0(\mathbb{P}^1,\mathcal{O}_{\mathbb{P}^1}(4)),\,wt\in H\}$.
Generically, one then has ${\rm dim}\,W'=2$ et $H=W'\cdot H^0(\mathbb{P}^1,\mathcal{O}_{\mathbb{P}^1}(2))$, and this concludes the proof of the claim.

 If the curve  $C$ is contained in  $Y$, the residual
curve of  $C$ in the intersection $Y\cap \Sigma$ is an elliptic curve of degree  $6$ in
$Y$.  These two  constructions are  inverse of each other, which provides a birational map  from the
space of pairs
 $\{(E,\Sigma),\,E\in M_{6,1},\,E\subset \Sigma\}$  to the space of pairs
$\{(C,\Sigma),\,C\in M_{6,0},\,C\subset \Sigma\}$, where $\Sigma\subset \mathbb{P}^4$ is a surface
of the type considered above, and the inclusions of the curves $C$, $E$ is $\Sigma$ are in the linear
systems described above.

 This construction thus provides rational curves  $R_C\subset M_{6,1}$ parameterized by  $C\in M_{6,0}$ and sweeping-out  $M_{6,1}$. To each curve  $C$ parameterized by  $M_{6,0}$, one associates the family $R_C$ of
elliptic curves  $E_t$ which are residual to $C$ in the intersection  $\Sigma_t\cap Y$, where $t$ runs over the $\mathbb{P}^1$ exhibited in the above claim.

Let $\Phi:M_{6,1}\dashrightarrow B$ be the maximal rationally connected fibration
of $M_{6,1}$.
We will use later on the
  curves $R_C$ introduced above  to prove the following Lemma \ref{surject51}.  Notice first
  that the divisor $D_{5,1}\subset M_{6,1}$
parameterizing singular elliptic curves  of degree $6$ contained in  $Y$, which are the union of an elliptic curve of
 degree  $5$ in $Y$ and of a line  $\Delta\subset Y$ meeting in one  point.  If $E$ is such a
 generically chosen  elliptic curve
  and    $S\in| \mathcal{O}_{Y}(2)|$ is a smooth  $K3$ surface containing  $E$, the linear system  $| \mathcal{O}_S(E)|$  contains  the line $\Delta$ in its base locus, and thus the
 construction of the associated vector bundle $\mathcal{E}$ fails. Consider however a curve
 $E'\subset S$ in the linear system
 $| \mathcal{O}_S(2-E)|$. Then  $E'$ is again a degree $6$ elliptic curve in
  $Y$ and $E'$  is smooth and nondegenerate for a generic choice of  $Y,E,S$ as above.
 The curves  $E'$ so obtained  are parameterized by a divisor
 $D'_{5,1}$ of $M_{6,1}$, which also has the following description : the curves $E'$ parameterized
  by $D'_{5,1}$ can also be characterized by the fact that
 they have  a trisecant line contained in  $Y$ (precisely the line  $\Delta$
 introduced above).

 \begin{lemm}\label{surject51} The restriction of $\Phi$ to $D'_{5,1}$ is surjective.
\end{lemm}
\begin{rema}{\rm It is not true that the map $\phi_{\mid D'_{5,1}}:D'_{5,1}\dashrightarrow M_9$ is dominating. Indeed, curves $E$ parameterized by $D'_{5,1}$ admit a trisecant line $L$ in $Y$. Let $\mathcal{E}$
be the associated vector bundle. Then because ${\rm det}\,\mathcal{E}=\mathcal{O}_Y(2)$ and $\mathcal{E}_{\mid L}$ has a nonzero section vanishing in three points, $\mathcal{E}_{\mid L}\cong\mathcal{O}_L(3)\oplus \mathcal{O}_L(-1)$. Hence the corresponding vector bundles $\mathcal{E}$ are not globally generated and thus are not generic.}
\end{rema}
Admitting Lemmas
\ref{lem8mars} and \ref{surject51}, the proof of Theorem \ref{propo61ratcon} is concluded in the following way:

 Lemma \ref{lem8mars} says that the rational map  $\Phi_{\mid D_{3,3}}$ is dominating. The dominating rational maps
$\Phi_{\mid D_{3,3}}$ and
$\Phi_{\mid D'_{5,1}}$, taking values in a non uniruled variety,
 factorize through the base of the MRC fibrations of
 $ D_{3,3}$ and  $D'_{5,1}$ respectively. The maximal rationally connected fibrations of the divisors
 $D'_{5,1}$ and
 $D_{5,1}$ have the same base. Indeed, the  construction of the  correspondence
  between $D_{5,1}$ and $D'_{5,1}$ is obtained by  liaison,
  and it follows that a  $\mathbb{P}^2$-bundle over  $D_{5,1}$ (parameterizing a reducible  elliptic curve $E$ of type
  $(5,1)$ and a rank  $2$ vector  space of
  quadratic polynomials  vanishing  on $E$) is
  birationally  equivalent to a  $\mathbb{P}^2$-fibration  over $D'_{5,1}$ defined in the same way. The MRC
  fibration of the divisor  $D_{5,1}$ is well understood by the  works
\cite{ilievmarku}, \cite{tikhomarku}. Indeed, these papers   show that
the morphism  $\phi_{5,1}$ induced by the Abel-Jacobi map of $Y$  on the family  $M_{5,1}$ of degree $5$ elliptic curves in  $Y$ is surjective with general fiber isomorphic to  $\mathbb{P}^5$: more precisely,
an elliptic curve  $E$ of degree  $5$ in  $Y$ is the zero locus of  a transverse section
 of a rank $2$ vector bundle  $\mathcal{E}$ on  $Y$, and
 the fiber of $\phi_{5,1}$ passing  through $[E]$ is for general $E$ isomorphic to   $\mathbb{P}(H^0(Y,\mathcal{E}))=\mathbb{P}^5$.
The vector bundle  $\mathcal{E}$, and the line  $L\subset Y$  being given, the set of elliptic curves
of type $(5,1)$ whose degree $5$ component is the zero locus of a section  of $\mathcal{E}$
 and the degree $1$ component is the line $L$ identifies to  $$\bigcup_{x\in L}\mathbb{P}(H^0(Y,\mathcal{E}\otimes\mathcal{I}_x)),$$ which is rationally connected.
 Hence $D_{5,1}$, and thus also  $D'_{5,1}$, is a  fibration over  $J(Y)\times F$
 with rationally connected general fiber, where
$F$ is the Fano surface of lines in  $Y$. As the variety  $J(Y)\times F$ is not uniruled, it must be the base of
 the maximal rationally connected fibration of  $D'_{5,1}$.

 From the above considerations, we conclude that $B$ is dominated by
 $J(Y)\times F$ and in particular  ${\rm dim}\,B\leq 7$.

Let us now consider the dominating rational map
$\Phi':=\Phi_{\mid D_{3,3}}:D_{3,3}\dashrightarrow B$.
Recall that  $D_{3,3}$ admits a surjective morphism $\chi$ to   $S^2\Theta$, the general fiber being isomorphic to
the second symmetric product  $S^2E_3$ of an elliptic curve of degree $3$. We postpone the proof of the following lemma:
\begin{lemm}\label{le8marsfacteur}
The rational map  $\Phi'$   factorizes through    $S^2\Theta$.
\end{lemm}

 This lemma implies that  $B$ is rationally dominated by $S^2\Theta$, hence a fortiori
by $\Theta\times \Theta$. We will denote  $\Phi'':\Theta\times\Theta\dashrightarrow B$
 the rational map obtained by  factorizing
 $\Phi'$ through $S^2\Theta$ and then by composing the resulting map  with the
 quotient map $\Theta\times \Theta\rightarrow S^2\Theta$.  We have a map
   $\phi_B:B\rightarrow J(Y)$ obtained as a factorization of  the map
  $\phi_{6,1}:M_{6,1}\rightarrow J(Y)$, using the fact that  the fibers of $\Phi:M_{6,1}\dashrightarrow B$ are rationally connected,
  hence contracted by $\phi_{6,1}$. Clearly, the composition of the
 map $\phi_B$ with  the map
 $\Phi''$ identifies to the sum map  $\sigma:\Theta\times \Theta\rightarrow J(Y)$.
 Its fiber over an element $a\in J^2(Y)$ thus identifies
 in a canonical way to the complete intersection  of the ample divisors
 $\Theta$ and  $ a-\Theta$ of $J(Y)$.   Consider the  dominating rational map
  fiberwise induced by  $\Phi''$ over $J(Y)$:
 $$\Phi''_a:\Theta\cap (a-\Theta)\dashrightarrow B_a,$$
 where  $B_a:=\phi_B^{-1}(a)$ has  dimension at most  $2$
 and is not uniruled. The proof is thus concluded with  Lemma \ref{thetaetsurface} below,
 which implies that  $B$ is in fact birationally isomorphic to  $J(Y)$. Indeed, this lemma tells that
 $\Phi''_a$ must be constant for very general  $Y$ and general  $a$, hence also
 for  general $Y$   and  $a$.
 \cqfd
 \begin{lemm} \label{thetaetsurface} If $Y$ is very general,  and $a\in J(Y)$ is general,  any
   dominating rational map $\Phi''_a$
 from $\Theta\cap (a-\Theta)$ to a smooth  non uniruled variety  $B_a$ of dimension
 $\leq 2$ is constant.
 \end{lemm}
 To complete the proof of Theorem \ref{propo61ratcon}, it remains to prove  Lemmas  \ref{lem8mars},
 \ref{surject51}, \ref{le8marsfacteur} and \ref{thetaetsurface}.

 \vspace{0.5cm}

{\bf Proof of Lemma \ref{lem8mars}.}  Let $E=C_3\cup C'_3$ be a general  elliptic curve $(3,3)$
contained in  $Y$. Then $E$ is contained in a smooth   $K3$ surface  in the linear system
$ | \mathcal{O}_Y(2)|$.  The linear  system $| \mathcal{O}_S(E)|$ has no base point in $S$, and thus the construction of the vector bundle $\mathcal{E}$ can be done as in the smooth general case, and furthermore it is easily verified that  $\mathcal{E}$ is stable on  $Y$. Hence  $\phi$ is
 well-defined at the point  $[E]$.
One verifies that  $H^0(Y,\mathcal{E})$ is of  dimension $4$ as in the general smooth case. As
${\rm dim}\,D_{3,3}=10$ and  ${\rm dim}\,M_9=9$, to show that  $\phi_{\mid D_{3,3}}$ is dominating, it suffices to show that the
general fiber of  $\phi_{\mid D_{3,3}}$ is of  dimension $\leq 1$. Assume to the contrary
 that this fiber, which is contained in  $\mathbb{P}(H^0(Y,\mathcal{E}))=\mathbb{P}^3$, is of dimension $\geq 2$. Recalling that $D'_{3,3}$ is not swept-out by rational  surfaces because it is fibered over
 $S^2\Theta$ into   surfaces
isomorphic to the second symmetric product of an elliptic curve,
  this fiber should be a surface of degree  $\geq 3$
in $\mathbb{P}(H^0(Y,\mathcal{E}))$.   A general line in
$\mathbb{P}(H^0(Y,\mathcal{E}))$ should then  meet this surface  in at least $3$ points counted with multiplicities. For such a line, take
the  $\mathbb{P}^1\subset \mathbb{P}(H^0(Y,\mathcal{E}))$ obtained by considering the base-point free
pencil
$| \mathcal{O}_S(E)|$. One verifies that this line is not  tangent to $D_{3,3}$ at the point $[E]$
and thus should meet  $D_{3,3}$ in another point. The surface $S$ should thus have a
 Picard number $\rho(S)\geq 4$, which is easily excluded by a dimension count for a general pair
$(E,S)$, where $E$ is an  elliptic curve of  type $(3,3)$ as above and
$S$  is the intersection of our general  cubic $Y$ and  a quadric in $\mathbb{P}^4$ containing
$E$.
\cqfd
For the proof of Lemma \ref{surject51} we will use the following elementary Lemma
 \ref{elementairesuivant} : Let
   $Z$ be a smooth projective variety, and $\Phi:Z\dashrightarrow B$
 its maximal rationally connected fibration. Assume given a family of rational curves $(C_t)_{t\in M}$
  sweeping-out $Z$.
 \begin{lemm}\label{elementairesuivant} Let $D\subset Z$ be a divisor such that, through a general point
  $d\in D$, there passes a curve $C_t,\,t\in M,$ meeting $D$ properly at $d$. Then
 $\Phi_{\mid D}:D\dashrightarrow B$ is dominating.

 \end{lemm}

  {\bf Proof of Lemma \ref{surject51}.} The proof is obtained by applying Lemma \ref{elementairesuivant} to
  $Z=M_{6,1}$ and to the previously constructed  family of rational
  curves  $R_C,\,C\in M_{6,0}$. It thus suffices to show that
   for general $Y$, if $[E]\in D'_{5,1}$ is a general point, there exists a (irreducible)  curve $R_C$   passing through $[E]$ and not contained in
   $D'_{5,1}$. Let us recall that $D'_{5,1}$
  is irreducible, and parameterizes
  elliptic curves of degree  $6$ in $Y$ admitting a trisecant line
  contained in $Y$. Starting from a curve
  $C$ of genus  $0$ and degree  $6$ which is nondegenerate in  $\mathbb{P}^4$, we constructed  a one
  parameter family $(\Sigma_t)_{t\in \mathbb{P}^1}$ of surfaces  of degree  $4$ containing $C$ and isomorphic to
   $\mathbb{P}^1\times \mathbb{P}^1$. For general  $t$, the trisecant lines of  $\Sigma_t$ in
  $\mathbb{P}^4$ form at most a $3$-dimensional variety and thus, for a general cubic  $Y\subset \mathbb{P}^4$, there exists
   not trisecant to $\Sigma_t$ contained in  $Y$. Let us choose a line
  $\Delta$  which is trisecant to  $\Sigma_0$
  but not  trisecant to $\Sigma_t$ for  $t$ close to $0$, $t\not=0$.
  A general cubic $Y$ containing $C$ and $\Delta$ contains then no trisecant to
  $\Sigma_t$ for $t$ close to $0$,  $t\not=0$.
  Let $E_0$ be the residual curve of  $C$ in  $\Sigma_0\cap Y$. Then $E_0$ has the line
  $\Delta\subset Y$ as a trisecant, and thus $[E_0]$ is a point of the divisor $D'_{5,1}$ associated to $Y$.
  However the rational curve $R_C\subset M_{6,1}(Y)$ parameterizing the
  residual curves   $E_t$ of $C$ in $\Sigma_t\cap Y,\,t\in \mathbb{P}^1$, is not contained in
   $D'_{5,1}$ because, by construction, no trisecant of $E_t$  is contained in $Y$,
  for $t$ close to $0$, $t\not=0$.
 \cqfd
 {\bf Proof of Lemma \ref{le8marsfacteur}.} Assume to the contrary that the dominating rational map
 $\Phi':D_{3,3}\dashrightarrow B$ does not factorize through  $S^2\Theta$. As the map
 $\chi:D_{3,3}\rightarrow S^2\Theta$ has for general fiber $S^2E_3$, where
 $E_3$ is an elliptic curve of degree $3$ in $\mathbb{P}^4$, and $B$ is not uniruled, $\Phi'$ factorizes
  through the corresponding fibration $\chi':Z\dashrightarrow S^2\Theta$, with  fiber $Pic^2E_3$. Let us denote $\Phi'' :Z\dashrightarrow B$ this second factorization. Recall that ${\rm dim}\,B\leq 7$ and that $B$ maps surjectively onto
  $J(Y)$ via the rational map $\phi_B:B\dashrightarrow J(Y)$ induced by  the Abel-Jacobi map
  $\phi_{6,1}: M_{6,1}\rightarrow J(Y)$.
  Clearly the elliptic curves $Pic^2E_3$ introduced above are contracted by
  the composite map $\phi_B\circ \Phi''$. If these curves are not contracted by $\Phi''$, the fibers
 $B_a:=\phi_B^{-1}(a),\,a\in J^2(Y),$ are swept-out by elliptic curves. As the fibers $B_a$ are
 either curves of genus $>0$
 or non uniruled surfaces  (because ${\rm dim}\,B\leq 7$ and $B$ is not uniruled), they
  have the property that there pass at most  finitely many  elliptic curves in a given deformation class through a general point
  $b\in B_a$. Under our assumption, there should thus
  exist :

  -  A variety $J'$, and
    a dominating morphism  $g:J'\rightarrow J(Y)$ such that the fiber  $J'_a$, for a general point $a\in J^2(Y)$,
     has dimension at most $1$
    and parameterizes elliptic curves in a given deformation class in the fiber $B_a$.

    - A
      family  $B'\rightarrow J'$ of elliptic curves, and dominating rational maps
     :
     $$\Psi: Z\dashrightarrow B',\,s:S^2\Theta\dashrightarrow J',\,r:Z'\dashrightarrow B$$

     such that  $Z$ is birationally equivalent to the fibered product
     $S^2\Theta\times_{J'}B'$, $r$ is generically finite, and $\Phi''=r\circ \Psi$.

  To get a contradiction out of this, observe that the family of elliptic curves $Z\rightarrow S^2\Theta$
  is also birationally the inverse image of a family of elliptic
 curves  parameterized by $G(2,4)$, via a natural rational map
 $f:S^2\Theta\dashrightarrow G(2,4)$. Indeed, an element of
 $S^2\Theta$ parameterizes  two linear systems  $| C_1|$, resp. $| C_2|$
 of rational curves of degree $3$  on two hyperplane sections
$S_1=H_1\cap Y$, resp. $S_2=H_2\cap Y$ of $Y$. The map  $f$ associates to the unordered pair
$\{| C_1|,| C_2|\}$ the plane $P:=H_1\cap H_2$. This plane $P\subset \mathbb{P}^4$ parameterizes the elliptic curve
 $E_3:=P\cap Y$, and there is a canonical isomorphism  $E_3\cong Pic^2E_3$ because
$deg\,E_3=3$.

The contradiction then comes from the following fact: for $j\in \mathbb{P}^1$, let us consider the divisor $D_j$ of $Z$
 which is swept-out by
elliptic curves of fixed modulus  determined by $j$, in the fibration
$Z\dashrightarrow S^2\Theta$. This divisor is the  inverse image by $f\circ\chi'$
of an ample divisor  on  $G(2,4)$ and furthermore it is also the  inverse image by
$\Psi$
of the similarly defined divisor of  $B'$.   As
     $r$ is generically finite and  $B$ is not uniruled, $B'$ is not uniruled.
     Hence the rational map $\Psi :Z\dashrightarrow B'$ is ``almost holomorphic'', that is
      well-defined along its general fiber.
One deduces from this that the divisor   $D_j$ has trivial restriction
 to a general fiber of $\Psi$. This fiber, which is of  dimension $\geq 2$ because
 ${\rm dim}\, Z=9$ and ${\rm dim}\,B'\leq 7$, must then be sent to a point in $G(2,4)$ by $f\circ \chi'$ and its image in
$S^2\Theta$ is thus contained in an union of components of fibers of $f$.
 Observing furthermore that  the general fiber of
${f}$  is irreducible of dimension $2$,
it follows that the rational dominating map
$s:S^2\Theta\dashrightarrow J'$ factorizes through
$f:S^2\Theta\dashrightarrow G(2,5)$. As $G(2,5)$ is rational, this contradicts the fact that
$J'$ dominates rationally $J(Y)$.

 \cqfd
 {\bf Proof of Lemma \ref{thetaetsurface}.}  Observe first that, as ${\rm dim}\, {\rm Sing}\,\Theta=0$,
 the general intersection $\Theta\cap (a-\Theta)$ is smooth.
 Assume first that ${\rm dim}\,B_a=1$. Then either $B_a$ is an elliptic curve or  $H^0(B_a,\Omega_{B_a})$ has dimension $\geq 2$.
 Notice that, by  Lefschetz hyperplane restriction  theorem, the Albanese variety
 of the complete intersection  $\Theta\cap (a-\Theta)$ identifies to $J(Y)$.  $B_a$ cannot be elliptic,
  otherwise $J(Y)$ would not be simple. But it is easy to prove by an infinitesimal variation of Hodge structure argument (cf. \cite[II, 6.2.1]{voisinbook}) that  $J(Y)$ is simple for very general  $Y$.
  If  ${\rm dim}\,H^0(B_a,\Omega_{B_a})\geq 2$, this provides two independent holomorphic $1$-forms  on
 $\Theta\cap (a-\Theta)$ whose  exterior product vanishes. These $1$-forms come, by
  Lefschetz hyperplane restriction theorem,
 from independent holomorphic  $1$-forms on  $J(Y)$, whose exterior product
 is a nonzero holomorphic  $2$-form on  $J(Y)$. The restriction of this $2$-form
 to  $\Theta\cap (a-\Theta)$ vanishes, which contradicts  Lefschetz hyperplane restriction
 theorem because
  ${\rm dim}\,\Theta\cap (a-\Theta)=3$. This case is thus excluded.

 Assume now ${\rm dim}\,B_a=2$. Then the surface $B_a$ has  nonnegative Kodaira dimension, because it is not uniruled. Let $L$ be the line bundle
  ${\Phi''_a}^*(K_{B_a})$ on $\Theta\cap (a-\Theta)$. The Iitaka dimension   $\kappa(L)$ of $L$ is nonnegative and there is a nonzero section of the bundle
 $\Omega_{\Theta\cap (a-\Theta)}^2(-L)$. Notice that by Grothendieck-Lefschetz theorem, $L$ is the restriction of a line bundle on
  $J(Y)$, hence must be numerically effective, because one can show, again  by an infinitesimal variation of Hodge structure argument, that ${\rm NS}(J(Y))=\mathbb{Z}$, for very general  $Y$.
 An easy argument using the conormal
  exact sequence of   $\Theta\cap (a-\Theta)$ in  $J^2(Y)$,
  vanishing theorems, and the fact that ${\rm dim}\,\Theta\cap (a-\Theta)\geq3$,
 allows to show that any such  section  is the restriction
of a section of $(\Omega_{J(Y)}^2)_{\mid \Theta\cap (a-\Theta)}(-L)$.
 As  $\kappa(L)\geq0$,  such a nonzero section
can exist only if   $L$ is trivial.
 Notice that, up to replacing  $(\Theta\cap (a-\Theta),\Phi''_a, B_a)$ by its   Stein factorization (or rather, the
  Stein factorization of a desingularized model  of $\Phi''_a$),
  one may assume that  $\Phi''_a$ has connected fibers. In this case, the triviality of $L$ implies the existence
 of a nonzero  section of  $K_{B_a}$ and the fact that  $h^0(B_a,K_{B_a})=1$.
 This is clear if $\Phi''_a$ is a morphism. Otherwise, choose a desingularization:
 $$\widetilde{\Phi''_a}:\widetilde{\Theta\cap (a-\Theta)}\rightarrow B_a$$
 of the rational map $\Phi''_a$. Then we know that $\widetilde{\Phi''_a}^*K_{B_a}$ has nonnegative Iitaka
 dimension and that it is trivial away from the exceptional divisor of the desingularization
 map
 $$\widetilde{\Theta\cap (a-\Theta)}\rightarrow \Theta\cap (a-\Theta).$$
 This is possible only if $\widetilde{\Phi''_a}^*K_{B_a}$ is effective supported on this exceptional divisor.
 But then it has only one section, which  provides a unique section of
 $K_{B_a}$, because we assumed that fibers of $\widetilde{\Phi''_a}$ are connected.

 Having this, we conclude that the  Hodge structure on
 $H^2(\Theta\cap (a-\Theta),\mathbb{Z})\cong H^2(J^2(Y),\mathbb{Z})$
 has    a Hodge substructure with  $h^{2,0}=1$. But this
 can be also excluded for  very general  $Y$
 by an infinitesimal variation of Hodge structure argument.

 \cqfd
This concludes the proof of  Theorem \ref{propo61ratcon}.

\subsection{Application: Integral Hodge classes on cubic threefolds fibrations over curves\label{seccub}}

 The main concrete application of Theorem \ref{propo61ratcon} concerns  fibrations
 $f:X\rightarrow \Gamma$ over a curve  $\Gamma$ with general fiber a cubic threefold
 in $\mathbb{P}^4$ (or the smooth projective models
of smooth cubic
hypersurfaces $X\subset\mathbb{P}^4_{\mathbb{C}(\Gamma)}$).
\begin{theo}\label{cubique} Let $f:X\rightarrow \Gamma$ be a cubic threefold  fibration over a smooth curve. If the
 fibers of $f$ have at worst ordinary quadratic singularities,
$Z^4(X)=0$ and  $H^3_{nr}(X,\mathbb{Z}/n\mathbb{Z})=0$  for any $n$.

\end{theo}
Note that the second statement is a  consequence of the first by Theorem \ref{ctVtheo}
 and the fact that  for $X$ as above, $CH_0(X)$ is supported on a curve.

Theorem \ref{cubique}  applies in particular to   cubic fourfolds. Indeed,  a smooth cubic hypersurface
$X$ in  $\mathbb{P}^5$ admits a   Lefschetz pencil of  hyperplane  sections. It thus becomes, after blowing-up the base-locus of this pencil,
birationally equivalent to a model of a cubic in  $\mathbb{P}^4_{\mathbb{C}(t)}$ which satisfies all
our hypotheses.  Theorem
\ref{cubique} thus provides $Z^4(X)=0$ for cubic fourfolds, a result first proved in   \cite{voisinjapjmath}.

\vspace{0.5cm}

{\bf Proof of Theorem  \ref{cubique}.} We first prove the theorem for
a fibration $X\rightarrow \Gamma$ whose general fiber satisfies  the conclusion of Theorem
 \ref{propo61ratcon}.

Cubic hypersurfaces in  $\mathbb{P}^4$ do not have  torsion in their degree  $3$
cohomology  by  Lefschetz  hyperplane restriction theorem. They furthermore  satisfy the vanishing conditions $H^i(X_t,\mathcal{O}_{X_t})=0$, $i>0$.

It thus suffices  to prove that  condition
\ref{itemref3} of Theorem \ref{theocritereabeljacobi} is satisfied.

The  fibers of $f$ being smooth hypersurfaces in $\mathbb{P}^4$,
$R^4f_*\mathbb{Z}$  is, over the regularity open set  $\Gamma_0$, the trivial
 local  system isomorphic to $\mathbb{Z}$. A section $\alpha$ of $R^4f_*\mathbb{Z}$ over $\Gamma_0$ is thus  characterized by  a number, its degree on the fibers  $X_t$ with respect to the   polarization $c_1(\mathcal{O}_{X_t}(1))$.
When the degree of $\alpha$ is $5$, one gets the desired family  $B_5$ and
 cycle $\mathcal{Z}_5$
using the results of
 \cite{ilievmarku}, in the following way: Iliev and  Markoushevich prove that if $Y$ is a smooth
  cubic threefold, denoting $M_{5,1}$ a desingularization  of the Hilbert
  scheme of
   degree  $5$, genus  $1$ curves in  $Y$,
  the map $\phi_{5,1}:M_{5,1}\rightarrow J(Y)$
  induced by the Abel-Jacobi map of  $Y$ is  surjective with rationally connected fibers.
 Taking for  $B_{5}$
 a desingularization   $\mathcal{M}_{5,1}$  of the  relative  Hilbert scheme of curves
   of degree  $5$ and genus  $1$ in the fibers of  $f$, and for  cycle
   $\mathcal{Z}_{5}$ the universal subscheme,
   the hypothesis \ref{itemref3} of Theorem   \ref{theocritereabeljacobi} is thus satisfied
for the degree $5$ section  of $R^4f_*\mathbb{Z}$.

 For the degree $6$ section $\alpha$ of $R^4f_*\mathbb{Z}$,  property  \ref{itemref3} is similarly a consequence of Theorem  \ref{propo61ratcon}. Indeed, it says that the general fiber of the morphism  $\phi_{6,1}: M_{6,1}\rightarrow J(Y)_6$ induced by the Abel-Jacobi map of  $Y$
on the family  $M_{6,1}$ of degree $6$ elliptic curves in $Y$ is rationally connected for general $Y$, which by assumption remains true for the general fiber of $f$.
We thus take as before for family  $B_{6}$
 a desingularization  $\mathcal{M}_{6,1}$ of the relative  Hilbert scheme of curves
   of degree  $6$ and genus  $1$ in the fibers of $f$, and for cycle
   $\mathcal{Z}_{6}$ the universal subscheme.

To conclude, let us show that property   \ref{itemref3} for the sections $\alpha_5$
of degree  $5$ and $\alpha_6$
of degree  $6$ imply property  \ref{itemref3} for any section  $\alpha$. To see this, let us
introduce
the cycle $h^2$ of codimension $2$ over  $X$, where $h\in Pic\,X$
induces by restriction to a general fibre $X_t$ the positive generator of  $Pic\, X_t$. The cycle $h^2$ is thus of degree $3$
on the fibers of $f$.
The degree of a  section $\alpha$ is congruent to  $5,\,-5$ or $6$ modulo $3$, and we can thus write
$\alpha=\pm  \alpha_5+ \mu\alpha_3$, or $\alpha=\alpha_6+\mu\alpha_3$ for an adequate number $\mu$.
 In the first case, consider the variety $B_{\alpha}=B_{5}$ and the cycle $$\mathcal{Z}_{\alpha}=\pm \mathcal{Z}_{5}+\mu (B_{5}\times_\Gamma  h^2)\subset B_{5}\times_\Gamma X.$$
In the second case, consider the  variety $B_{\alpha}=B_{6}$ and the cycle $$\mathcal{Z}_{\alpha}= \mathcal{Z}_{6}+\mu (B_{6}\times\Gamma  h^2)\subset B_{6}\times_\Gamma X.$$
 It is clear that the pair
 $(B_{\alpha},\mathcal{Z}_{\alpha})$ satisfies the condition
 \ref{itemref3} of Theorem  \ref{theocritereabeljacobi}.

The general case, where the general fiber of $f:X\rightarrow \Gamma$ is not assumed to satisfy the conclusion
of Theorem \ref{propo61ratcon}, is treated in the same way, replacing the reference to Theorem
\ref{theocritereabeljacobi} by a reference to its variant (Theorem \ref{variantcritere}) and applying
the  arguments  above, reducing the verification of the assumption (ii)
  in Theorem \ref{variantcritere} to the case of the sections $\alpha$ of degree $5$ or $6$.
\cqfd

\section{Structure of the Abel-Jacobi map and  decomposition of the diagonal\label{vraiesecdecomp}}

\subsection{Relation between Questions \ref{question1} and
 \ref{question2}}
We establish in this subsection the following  relation between Questions \ref{question1} and
 \ref{question2}:
 \begin{theo} \label{theonouveau} Assume that Question \ref{question2} has a positive answer for $Y$ and that
the intermediate Jacobian of
$Y$ admits a  $1$-cycle $\Gamma$ of class $\frac{[\Theta]^{g-1}}{(g-1)!}$, $g=dim\,J(Y)$. Then
Question \ref{question1} also has an affirmative answer for $Y$.
\end{theo}
{\bf Proof.} There exists by assumption a variety $B$, and a codimension $2$ cycle
$Z\subset B\times Y$ cohomologous to $0$ on fibers $b\times Y$, such that the morphism
$$\phi_Z:B\rightarrow J(Y)$$
induced by the Abel-Jacobi map of $Y$ is surjective with rationally connected general fibers.
Consider the $1$-cycle $\Gamma$ of $J(Y)$. We may assume by moving lemma, up to changing the representative of
$\Gamma$ modulo homological equivalence, that $\Gamma=\sum_in_i\Gamma_i$
 where, for each component $\Gamma_i$ of
the support of $\Gamma$, the general fiber of $\phi_Z$ over $\Gamma_i$ is rationally connected.
According to \cite{GHS}, the inclusion $j_i:\Gamma_i\hookrightarrow J(Y)$ has then a lift
$\sigma_i:\Gamma_i\rightarrow B$. Denote  $Z_i\subset \Gamma_i\times Y$ the codimension $2$-cycle
$(\sigma_i,Id_Y)^*Z$.
Then the morphism $\phi_i:\Gamma_i\rightarrow J(Y)$ induced by the Abel-Jacobi map is
equal to $j_i$.

Observe now the following: there is the Pontryagin product $*$ on cycles of
$J(Y)$. Another way of rephrasing the condition that
the class of $\Gamma$ is equal to $\frac{[\Theta]^{g-1}}{(g-1)!}$, $g=dim\,J(Y)$, for some
principal polarization $\Theta$  is to say that
$$(\Gamma)^{*g}=J(Y)$$
as dimension $g$ cycles of $J(Y)$.

For each $g$-uple of components $(\Gamma_{i_1},\ldots,\Gamma_{i_g})$ of ${\rm Supp}\,\Gamma$,
consider $\Gamma_{i_1}\times\ldots\times \Gamma_{i_g}$, and the codimension $2$-cycle
$$Z_{i_1,\ldots,i_g}:=(pr_1,Id_Y)^*Z_{i_1}+\ldots+(pr_g,Id_Y)^*Z_{i_g}\subset \Gamma_1\times\ldots\times \Gamma_g\times Y.$$
The codimension $2$-cycle
$$\mathcal{Z}:=\sum_{i_1,\ldots,i_g}n_{i_1}\ldots n_{i_g}Z_{i_1,\ldots,i_g}\subset (\sqcup\Gamma_i)^g\times Y$$
is invariant under the symmetric group $\mathfrak{S}_g$ acting on the  factor
$(\sqcup\Gamma_i)^g$ in the product $(\sqcup\Gamma_i)^g\times Y$. The part of $\mathcal{Z}$ dominating
 over a component of $(\sqcup\Gamma_i)^g$, (which is the only one we are interested in)  is then the pull-back of  a codimension $2$-cycle $\mathcal{Z}_{sym}$
on
$(\sqcup\Gamma_i)^{(g)}\times Y$. Consider now the sum map
$$\sigma:(\sqcup\Gamma_i)^{(g)}\rightarrow J(Y).$$
Let $\mathcal{Z}_J:=(\sigma,Id)_*\mathcal{Z}_{sym}\subset J(Y)\times Y$.
The proof concludes with the following:
\begin{lemm} The Abel-Jacobi map :
$$\phi_{\mathcal{Z}_J}:J(Y)\rightarrow J(Y)$$
is equal to $Id_{J(Y)}$.
\end{lemm}
{\bf Proof.} Instead of the symmetric product $(\sqcup\Gamma_i)^{(g)}$ and the
descended cycle $\mathcal{Z}_{sym}$, consider the
product $(\sqcup\Gamma_i)^{g}$, the cycle $\mathcal{Z}$ and the sum map
$$\sigma':(\sqcup\Gamma_i)^{g}\rightarrow J(Y).$$
Then we have $(\sigma',Id)_*\mathcal{Z}=g!(\sigma,Id)_*\mathcal{Z}_{sym}$ in $CH^2(J(Y)\times Y)$, so that
writing $\mathcal{Z}'_J:=(\sigma',Id)_*\mathcal{Z}$, it suffices to prove that:
\begin{eqnarray}\label{truc7avril}g!\,Id_{J(Y)}=\phi_{\mathcal{Z}'_J}:J(Y)\rightarrow J(Y).
\end{eqnarray}

This is done as follows: let $j\in J(Y)$ be a general point, and let $\{x_1,\ldots,x_N\}$ be the fiber of
$\sigma'$ over $j$. Thus each $x_l$ parameterizes a $g$-uple
$(i^l_1,\ldots,i^l_g)$ of components of ${\rm Supp}\,\Gamma$, and points
$\gamma^l_{i_1},\ldots,\gamma^l_{i_g}$ of $\Gamma_{i_1},\ldots,\Gamma_{i_g}$ respectively,
such that:
\begin{eqnarray}\label{eqn7avril}\sum_{1\leq k\leq g}\gamma^l_{i_k}=j.
\end{eqnarray}
On the other hand, recall that
\begin{eqnarray}\label{eqnbis7avril}\gamma^l_{i_k}=AJ_Y(Z_{{i^l_k},\gamma^l_{i_k}}).
\end{eqnarray}
It follows from (\ref{eqn7avril}) and (\ref{eqnbis7avril}) that for each $l$, we have:
\begin{eqnarray}\label{eqnter7avril}
AJ_Y(\sum_{1\leq k\leq g}Z_{i_k}^l,\gamma^l_{i_k} ) =AJ_Y(Z_{i^l_1,\ldots,i^l_g,(\gamma^l_{i_1},\ldots,\gamma^l_{i_g})})=j.
\end{eqnarray}
Recall now that
$\Gamma=\sum_in_i\Gamma_i$ and recall that
$(\Gamma)^{*g}=J(Y)$, which can be translated by saying that
$$\sigma'_*(\sum_{i_{1},\ldots,i_g}n_{i_1}\ldots n_{i_g} \Gamma_{i_1}\times\ldots\times\Gamma_{i_g})=
g! J(Y).$$
This exactly says that
$$\sum_{1\leq l\leq N} \sum_{i^l_{1},\ldots,i^l_g}n_{i^l_1}\ldots n_{i^l_g}=g!$$
which together with (\ref{eqnter7avril}) proves (\ref{truc7avril}).
\cqfd
The proof of Theorem \ref{theonouveau} is now complete.
\cqfd
\begin{rema}{\rm  The question whether the intermediate Jacobian of
$Y$ admits a  $1$-cycle $\Gamma$ of class $\frac{[\Theta]^{g-1}}{(g-1)!}$, $g=dim\,J(Y)$
is unknown even for the cubic threefold. However it  has a positive answer for $g\leq 3$ because any principally polarized abelian variety of dimension $\leq3$ is the Jacobian of a curve.
}
\end{rema}
\subsection{Decomposition of the diagonal modulo homological equivalence \label{secdecomp}}
 This section is devoted to the study of
 assumption \ref{itemref3}   in Theorem  \ref{theocritereabeljacobi}, (see also Question \ref{question2}). We study here this condition when the family  $f:X\rightarrow \Gamma$ is a product $Y\times\Gamma$, where
$Y$ is a smooth projective threefold staisfying $Z^4(Y)=0$
(for example, by theorem \ref{theovoisinuniruled}, if $Y$ is a uniruled threefold).
First of all, we observe that assumption  \ref{itemref3} in theorem  \ref{theocritereabeljacobi} is equivalent to the  following condition on the   $3$-fold
$Y$ with $H^2(Y,\mathcal{O}_Y)=0$ and $Z^4(Y)=0$.

\vspace{0.5cm}

(*) {\it There exists a smooth variety $B$ and a family of  $1$-cycles
$\mathcal{Z}\subset B\times Y$ homologous to  $0$ in $Y$, such that the Abel-Jacobi map $AJ_Y$
induces  a surjective  morphism $\phi_{\mathcal{Z}}:B\rightarrow J^2(Y)$ with rationally connected
general fiber.}

\vspace{0.5cm}

Indeed, as $Y$ satisfies
$H^2(Y,\mathcal{O}_Y)=0$ and $Z^4(Y)=0$,
any  class $\alpha\in H^4(Y,\mathbb{Z})$ is algebraic, that is, there exists $z_\alpha\in CH^2(Y)$, such that $[z]=\alpha$. If  $B$  and
$\mathcal{Z}$ exist as above, the family  $B\times\Gamma $ and the cycle $\mathcal{Z}\times \Gamma+B\times z_\alpha\times\Gamma\subset B\times Y\times \Gamma$ induce a surjective
map with rationally connected general fiber
$$\phi_\alpha:B\times\Gamma\rightarrow J(Y)_\alpha\times\Gamma=\mathcal{J}_\alpha$$ and thus
property \ref{itemref3} is satisfied by $Y\times \Gamma$.

Assume furthermore  $Y$ satisfies $CH_0(Y)=\mathbb{Z}$.
 Bloch-Srinivas decomposition of the diagonal \cite{blochsrinivas} says that there exists
a nonzero integer $N$ such that, denoting $\Delta_Y\subset Y\times Y$ the diagonal :
\begin{eqnarray}\label{decompdiag}N\Delta_Y=Y\times y+Z\,\,{\rm in}\,\, CH^3(Y\times Y),
\end{eqnarray}
where $y$ is any  point of  $Y$ and the  support of $Z$ is contained in $ D\times Y$, $D\varsubsetneqq Y$.

   We wish to study the invariant of $Y$ defined as the gcd of the integers $N$ appearing above. This is a birational invariant of  $Y$. One can also consider
the decomposition (\ref{decompdiag}) modulo  homological equivalence, and our results relate
the  triviality  of this invariant, that is the existence of an integral homological decomposition
 of the diagonal as in (\ref{decompdiag}), to  condition (*) (among other things).
We first recall that the intermediate  Jacobian
$$J(Y)=H^3(Y,\mathbb{C})/(F^2H^3(Y,\mathbb{C})\oplus H^3(Y,\mathbb{Z})/{\rm torsion}$$ of a smooth projective
threefold $Y$ with $H^3(Y,\mathcal{O}_Y)=0$ is an abelian variety, which is principally
polarized by the unimodular intersection form on the lattice
$H^3(Y,\mathbb{Z})/{\rm torsion}$. We denote by $\Theta$ a divisor  in the numerical divisor class defined by this
 principal polarization.

\begin{theo} \label{coup1}Let $Y$ be a smooth projective  $3$-fold.
Assume  $Y$ admits an integral homological decomposition of the  diagonal as in   (\ref{decompdiag}); then we have:

(i)  $H^*(Y,\mathbb{Z})$ is without  torsion.

(ii)
 $Z^4(Y)=0$.

(iii)  $H^i(Y,\mathcal{O}_Y)=0,\,\forall i>0$ and
    $Y$ satisfies  condition (*).
 \end{theo}

{\bf Proof.}  There  exists  by assumption  $D\varsubsetneqq Y$, which one may assume
 of pure   dimension  $2$, and  $Z\in CH^3(Y\times Y)$ such that
the  support of $Z$ is contained in  $ D\times Y$, and
\begin{eqnarray}\label{decompdiaghomo}[\Delta_Y]=[Y\times y]+[Z]\,\,{\rm in}\,\, H^6(Y\times Y,\mathbb{Z}).
\end{eqnarray}

Codimension $3$ cycles $z$ of    $Y\times Y$ act on  $H^*(Y,\mathbb{Z})$
and on the intermediate Jacobian of  $Y$, and this
action, which we will denote
$$ z^*:H^*(Y,\mathbb{Z})\rightarrow H^*(Y,\mathbb{Z}),\, z^*:J(Y)\rightarrow J(Y)$$
depends only on the cohomology class of $z$.
As the diagonal  of  $Y$ acts by the identity map on  $H^*(Y,\mathbb{Z})$ for $*>0$ and on  $J(Y)$, and the cycle
 $Y\times y$ acts trivially on both of these groups, one  concludes that
\begin{eqnarray}\label{eqn16mars1} Id_{H^*(Y,\mathbb{Z})}=Z^*:H^*(Y,\mathbb{Z})\rightarrow H^*(Y,\mathbb{Z}),\,{\rm for}\,\,*>0,
\end{eqnarray}
\begin{eqnarray}\label{eqn15mars1} Id_{J(Y)}=Z^*:J(Y)\rightarrow J(Y).
\end{eqnarray}
 Let $\tau:\widetilde{D}\rightarrow Y$ be a desingularization of $D$ and ${i}_{\widetilde{D}}=i_D\circ\tau:\widetilde{D}\rightarrow Y$.
The part of the cycle  $Z$ which dominates $D$ can be lifted to a cycle $\widetilde{Z}$
in   $\widetilde{D}\times Y$, and the remaining part acts trivially on $H^*(Y,\mathbb{Z})$ for $*\leq3$ for codimension reasons. Thus
the map $Z^*$ acting on $H^*(Y,\mathbb{Z})$ for $*\leq3$ can  be written as
\begin{eqnarray}\label{eqn15mars2}Z^*={i}_{\widetilde{D}*}\circ \widetilde{Z}^*.
\end{eqnarray}
We note that the action of  $\widetilde{Z}^*$ on  cohomology sends $H^*(Y,\mathbb{Z}),\,*\leq3,$ to
$H^{*-2}(\widetilde{D},\mathbb{Z}),\,*\leq 3$. The last groups have no  torsion. It follows that
$\widetilde{Z}^*$ annihilates the torsion of $H^*(Y,\mathbb{Z}),\,*\leq3$. Formula (\ref{eqn16mars1})
implies then that  $H^*(Y,\mathbb{Z})$ has no torsion for $*\leq3$.

To deal with the torsion of $H^*(Y,\mathbb{Z})$ with $*\geq4$, we rather use the action
$z_*$ of $z$ on $H^*(Y,\mathbb{Z}),\,*=4,\,5$. This action again factors through
$\widetilde{Z}_*$, (indeed, the part of $Z$ which does not dominate $D$ acts trivially on these groups for dimension reasons,) hence through the restriction map:
$$H^*(Y,\mathbb{Z})\rightarrow H^*(\widetilde{D},\mathbb{Z}).$$
As ${\rm dim}\widetilde{D}\leq 2$, $H^*(\widetilde{D},\mathbb{Z})$ has no torsion for $*=4,\,5$.
It follows that we also have $H^*(Y,\mathbb{Z})_{tors}=0$ for $*=4,\,5$.
This proves (i).

To prove (ii), let us consider again the action $Z_*=Id$ on the cohomology $H^4(Y,\mathbb{Z})$.
Observe again that the part $Z'$ of $Z$ not dominating $Z$  has a trivial action on $H^4(Y,\mathbb{Z})$, while the dominating part lifts as above to
a cycle $\widetilde{Z}$
in   $\widetilde{D}\times Y$. Then we find that
$Id_{H^4(Y,\mathbb{Z})}$ factors as  $\widetilde{Z}_*\circ {i}_{\widetilde{D}}^*$, hence through the restriction map ${i}_{\widetilde{D}}^*:H^4(Y,\mathbb{Z})\rightarrow H^4(\widetilde{D},\mathbb{Z})$.
As ${\rm dim}\,\widetilde{D}=2$, the group on the right is generated by classes of algebraic cycles,
and thus $H^4(Y,\mathbb{Z})$ is generated by classes of algebraic cycles on $Y$. In particular
$Z^4(Y)=0$.

(iii) The vanishing $H^i(Y,\mathcal{O}_Y)=0,\,\forall i>0$ is a well-known consequence of the
decomposition of the diagonal (\ref{decompdiag}) (cf. \cite{blochsrinivas}, \cite[II,10.2.2]{voisinbook}). This is important to guarantee that $J(Y)$ is an abelian variety.

It is well-known, and this is a consequence of  Lefschetz theorem on   $(1,1)$-classes, that there exists
a universal divisor $\mathcal{D}\subset {\rm Pic}^0(\widetilde{D})\times \widetilde{D}$, such that
 the induced morphism : $\phi_{\mathcal{D}}:{\rm Pic}^0(\widetilde{D})\rightarrow {\rm Pic}^0(\widetilde{D})$ is the identity. On the other hand, we have the morphism
 $$\widetilde{Z}^*:J(Y)\rightarrow J^1(\widetilde{D})=Pic^0(\widetilde{D}),$$
 which is a morphism of abelian varieties.

Let us consider the cycle $$\mathcal{Z}:=(Id_{J(Y)},{i}_{\widetilde{D}})_*\circ(\widetilde{Z}^*, Id_{\widetilde{D}})^*(\mathcal{D})\in CH^2(J(Y)\times Y).$$ Then $\phi_{\mathcal{Z}}: J(Y)\rightarrow J(Y)$ is equal to
$${i}_{\widetilde{D}*}\circ \phi_{\mathcal{D}}\circ \widetilde{Z}^*: J(Y)\rightarrow J^1(\widetilde{D})\rightarrow J^1(\widetilde{D})\rightarrow J(Y).$$
 As $\phi_{\mathcal{D}}$ is the identity map acting on   ${\rm Pic}^0(\widetilde{D})$ and ${i}_{\widetilde{D}*}\circ \widetilde{Z}^*$ is
the identity map acting on  $J(Y)$ according to (\ref{eqn15mars1}) and (\ref{eqn15mars2}), one concludes that
$\phi_{\mathcal{Z}}=Id: J(Y)\rightarrow J(Y)$, which of course implies (*).
\cqfd

\begin{rema}{\rm The proof even shows that under our assumptions, Question \ref{question1} has an affirmative answer for $Y$.}
\end{rema}
Recall from section \ref{subsectionintro} that for any smooth projective complex variety $Y$, $Z^4(Y)$ is a quotient of
$H^3_{nr}(Y,\mathbb{Q}/\mathbb{Z})$, if $Z^4(Y)$ is of torsion (which is always the case
in dimension $3$).
We can now get a  result better than (ii) if instead of considering the decomposition of the diagonal modulo homological equivalence, we consider it modulo algebraic equivalence. For example, we have the following statement, which
improves (iii) above:
\begin{prop}\label{coup3} Let $Y$ be a smooth projective  $3$-fold.
Assume  $Y$ admits an integral  decomposition of the  diagonal as in   (\ref{decompdiag}), modulo  algebraic
equivalence. Then $H^3_{nr}(Y,A)=0$ for any abelian group $A$.
\end{prop}
{\bf Proof.} We use the fact that correspondences modulo algebraic equivalence act on
cohomology groups $H^p(X_{Zar},\mathcal{H}^q)$. We refer to the appendix of \cite{ctvoisin} for this fact which is precisely stated as follows:
\begin{prop} If $X$, $Y$ are smooth projective and $Z\subset X\times Y$ is a cycle defined up to algebraic equivalence, satisfying
${\rm dim}\,Y-dim\,Z=r$, then there is an induced morphism
$$Z_*:H^p(X_{Zar},\mathcal{H}^q(A))\rightarrow H^{p+r}(Y_{Zar},\mathcal{H}^{q+r}(A)).$$
These actions are compatible with the composition of correspondences.
\end{prop}
Assume now that $Y$ admits a decomposition of the  diagonal of the form
$$\Delta_Y=Z_1+Z_2\,\,{\rm in}\,\, CH(Y\times Y)/{\rm alg},$$
where $Z_1\subset D\times Y$ for some $D\subsetneqq Y$ and $Z_2\subset Y\times y_0$.
The diagonal acts on $H^3_{nr}(Y,A)$ as the identity. Thus we get
$$Id_{H^3_{nr}(Y,A)} = Z_{1*}+Z_{2*}:H^3_{nr}(Y,A)\rightarrow H^3_{nr}(Y,A).$$
We observe now (by introducing again a desingularization $\widetilde{D}$ of $D$) that
$Z_{1*}=0$ on $H^3_{nr}(Y,A)$, because $Z_{1*}$ factors through the
restriction map
$$H^3_{nr}(Y,A)\rightarrow H^3_{nr}(\widetilde{D},A)$$
and the group on the right is $0$, because ${\rm dim}\,\widetilde{D}\leq 2$. Furthermore
$Z_{2*}$ also vanishes on $H^3_{nr}(Y,A)=H^0(Y_{Zar},\mathcal{H}^3(A))$ because $Im\,Z_{2*}$ factors through
$j_{y_0*}$ which shifts $H^p(\mathcal{H}^q),\,p\geq0,\,q\geq0,$ to $H^{p+3}(\mathcal{H}^{q+3})$.

Thus $Id_{H^3_{nr}(Y,A)}=0$ and ${H^3_{nr}(Y,A)}=0$.
\cqfd
 A partial converse to Theorem \ref{coup1} is as follows.

   \begin{theo} \label{coupe2} Assume that   $Y$ satisfies

   i) $CH_0(Y)=\mathbb{Z}$.

    ii) $Z^4(Y)=0$.

    iii)  $H^*(Y,\mathbb{Z})$ has no torsion and the intermediate Jacobian of
$Y$ admits a $1$-cycle $\Gamma$ of class $\frac{[\Theta]^{g-1}}{(g-1)!}$, $g=dim\,J(Y)$.

Then  condition (*) implies the
existence of an integral homological  decomposition of the diagonal as in (\ref{decompdiag}).
\end{theo}
\begin{rema} {\rm The condition $Z^4(Y)=0$ is satisfied if $Y$ is uniruled. Condition ii) is satisfied if $Y$ is rationally connected. For a rationally connected threefold $Y$, $H^*(Y,\mathbb{Z})$ is without torsion
if $H^3(Y,\mathbb{Z})$ is without torsion.}
\end{rema}
\begin{rema} {\rm That the integral decomposition of the diagonal as in  (\ref{decompdiag}) and
 in the Chow group $CH(Y\times Y)$ implies that $H^3(Y,\mathbb{Z})$ has no torsion was observed by Colliot-Th\'{e}l\`{e}ne.
}
\end{rema}
{\bf Proof of Theorem \ref{coupe2}.}  When the integral cohomology of  $Y$ has no torsion, the class
of the diagonal  $[\Delta_Y]\in H^6(Y\times Y,\mathbb{Z})$
has an integral  K\"unneth decomposition.
$$[\Delta_Y]=\delta_{6,0}+\delta_{5,1}+\delta_{4,2}+\delta_{3,3}+\delta_{2,4}+\delta_{1,5}+\delta_{0,6},$$
where $\delta_{i,j}\in H^i(Y,\mathbb{Z})\otimes H^j(Y,\mathbb{Z})$. The class
$\delta_{0,6}$ is the  class of  $Y\times y$ for any  point $y$ of $Y$.
 As we assumed $CH_0(Y)=0$, we have
\begin{eqnarray}\label{derder}H^1(Y,\mathcal{O}_Y)= 0,\,H^2(Y,\mathcal{O}_Y)=0.
 \end{eqnarray}
 The first condition implies that the
groups  $H^1(Y,\mathbb{Q})$ and $H^5(Y,\mathbb{Q})$ are trivial, hence the groups
$H^1(Y,\mathbb{Z})$ and $H^5(Y,\mathbb{Z})$ must be trivial since they have no torsion by assumption.
It follows that  $\delta_{5,1}=\delta_{1,5}=0$.

Next, the  second vanishing condition in (\ref{derder}) implies that the Hodge structure on
$H^2(Y,\mathbb{Q})$, hence also on $H^4(Y,\mathbb{Q})$ by duality, is trivial. Hence
 $H^4(Y,\mathbb{Z})$ and $H^2(Y,\mathbb{Z})$ are generated  by Hodge classes, and because we assumed $Z^4(Y)=0$,
 it follows that $H^4(Y,\mathbb{Z})$ and $H^2(Y,\mathbb{Z})$ are generated  by
cycle classes. From this one concludes that
$\delta_{4,2}$ et $\delta_{2,4}$ are represented by algebraic cycles  whose support
does not dominate  $Y$ by the first projection. The same is true for  $\delta_{6,0}$ which is the  class
 of $y\times Y$. The existence of a  decomposition as in  (\ref{decompdiaghomo}) is thus equivalent to the fact that there   exists
a cycle $Z\subset Y\times Y$, such that
  the support of  $Z$ is contained in  $ D\times Y$, with $D\varsubsetneqq Y$, and
  $Z^*:H^3(Y,\mathbb{Z})\rightarrow H^3(Y,\mathbb{Z})$ is the identity map. This last
   condition is indeed equivalent to the fact
  that the component of  type $(3,3)$ of  $[Z]$ is equal to $\delta_{3,3}$.

Let now  $\Gamma=\sum_in_i\Gamma_i$ be a  $1$-cycle of  $J(Y)$ of class
$\frac{[\Theta]^{g-1}}{(g-1)!}$,  where   $\Gamma_i\subset J(Y)$ are smooth curves
in general position. If  (*) is satisfied, there exist a variety $B$ and a family  $\mathcal{Z}\subset B\times Y$
of $1$-cycles homologous to  $0$ in $Y$, such that $\phi_{\mathcal{Z}}: B\rightarrow J(Y)$ is surjective with rationally connected general fiber. By  \cite{GHS}, $\phi_{\mathcal{Z}}$ admits a section
$\sigma_i$
 over each $\Gamma_i$, and $(\sigma_i,Id)^*\mathcal{Z}$ provides a family $Z_i\subset \Gamma_i\times Y$
of $1$-cycles
homologous to  $0$ in $Y$, parameterized by  $\Gamma_i$, such that
$\phi_{Z_i}:\Gamma_i\rightarrow J(Y)$ identifies to the  inclusion of
$\Gamma_i$ in  $J(Y)$.

Let us consider the cycle $Z\in CH^3(Y\times Y)$ defined by
$$Z=\sum_i n_iZ_i\circ {^tZ_i}.$$
The proof that the cycle $Z$ satisfies the desired property  is then given in the
following  Lemma \ref{lecalcul}.

\cqfd
\begin{lemm}\label{lecalcul} The map  $Z^*:H^3(Y,\mathbb{Z})\rightarrow H^3(Y,\mathbb{Z})$ is the identity map.
\end{lemm}
{\bf Proof.} We have
$$Z^*=\sum_in_i{^tZ_i}^*\circ Z_i^*.$$
Let us study the composite map
$$^tZ_i^*\circ Z_i^*: H^3(Y,\mathbb{Z})\rightarrow H^1(\Gamma_i,\mathbb{Z})\rightarrow H^3(Y,\mathbb{Z}).$$
Recalling that $Z_i\in CH^2(\Gamma_i\times Y)$ is the restriction to $\sigma(\Gamma_i)$ of
$\mathcal{Z}\in CH^2(\Gamma_i\times Y)$, one finds that
this composite map can also be written as:
$$^tZ_i^*\circ Z_i^*= {
^t\mathcal{Z}^*}\circ \cup[\sigma(\Gamma_i)]\circ \mathcal{Z}^*.$$
One uses for this  the fact that,  denoting  $j_i:\sigma(\Gamma_i)\rightarrow B$ the inclusion,
the composition
$$j_{i*}\circ j_i^*: H^1(B,\mathbb{Z})\rightarrow H^{2n-1}(B,\mathbb{Z}),\,n=dim\,B$$ is equal to
$\cup[\sigma(\Gamma_i)]$.

We thus obtain:
$$Z^*=
{^t\mathcal{Z}^*}\circ \cup(\sum_in_i[\sigma(\Gamma_i)])\circ \mathcal{Z}^*.$$
But we know that the map  $\phi_\mathcal{Z}: B\rightarrow J(Y)$
has rationally connected fibers,
and thus induces isomorphisms:
\begin{eqnarray}\label{chou}\phi_\mathcal{Z}^*:H^1(J(Y),\mathbb{Z})\cong H^1(B,\mathbb{Z}),\,\phi_{\mathcal{Z}*}: H^{2n-1}(B,\mathbb{Z})\cong H^{2g-1}(J(Y),\mathbb{Z}).
\end{eqnarray}
Via these  isomorphisms,
$\mathcal{Z}^*$ becomes the canonical isomorphism  $$H^3(Y,\mathbb{Z})\cong H^1(J(Y),\mathbb{Z})$$
(which  uses the fact that $H^3(Y,\mathbb{Z})$ is torsion free) and $^t\mathcal{Z}^*$ becomes
the dual canonical isomorphism
$$H^{2g-1}(J(Y),\mathbb{Z})\cong H^3(Y,\mathbb{Z}).$$
Finally, the map $\cup(\sum_in_i[\sigma(\Gamma_i)]):H^1(B,\mathbb{Z})\rightarrow H^{2n-1}(B,\mathbb{Z})$
identifies via (\ref{chou}) to the map
$\cup(\sum_in_i[(\Gamma_i)]):H^1(J(Y),\mathbb{Z})\rightarrow H^{2g-1}(J(Y),\mathbb{Z})$.
Recalling that $\sum_in_i[(\Gamma_i)]=\frac{[\Theta]^{g-1}}{(g-1)!}$, we identified
$Z^*:H^3(Y,\mathbb{Z})\rightarrow H^3(Y,\mathbb{Z})$ to the composite map
$$H^3(Y,\mathbb{Z})\cong H^1(J(Y),\mathbb{Z})\stackrel{\cup \frac{[\Theta]^{g-1}}{(g-1)!}}{\rightarrow}
H^{2g-1}(J(Y),\mathbb{Z})\cong H^3(Y,\mathbb{Z}),$$
where the last isomorphism  is Poincar\'{e} dual   of the first,
and using the definition of the  polarization $\Theta$ on $J(Y)$ (as being given
  by Poincar\'{e} duality on $Y$: $H^3(Y,\mathbb{Z})\cong H^3(Y,\mathbb{Z})^*$) we find that this composite map is the identity.
\cqfd

\end{document}